\documentclass[aap,MSNbibl,seceqn,dvips]{arximspdf}
\usepackage{graphicx,url,breakurl}
\doi{10.1214/13-AAP959} 
\volume{24}
\issue{4}
\pubyear{2014}
\firstpage{1652}
\lastpage{1690}

\makeatletter
\newcommand{\rrvert}{\vert}
\newcommand{\llvert}{\vert}
\def\cal{\mathcal}
\def\emphh{}
\newcommand{\eqref}[1]{(\ref{#1})}
\newtheorem{theorem}{Theorem}[section]
\newtheorem{lemma}{Lemma}[section]
\newproclaim{example}{Example}[section]
\newproclaim{definition}{Definition}[section]
\newproclaim{remark}{Remark}[section]
\newtheorem{proposition}{Proposition}[section]
\newtheorem{corollary}{Corollary}[section]

\def\st{\Sigma(T-t)}
\def\E{\mathbb{E}}
\def\F{\mathcal{F}}
\def\P{\mathbb{P}}
\def\R{\mathbb{R}}
\def\S{\mathcal{S}}
\def\C{\mathcal{C}}
\def\D{\mathcal{D}}
\def\A{\mathcal{A}}
\def\somt{\sum_{\tau_{i-1}^n< t}}
\def\supi{\sup_{1\leq i\leq N^n_{T}}}
\def\supit{\sup_{\tau^n_{i-1} \leq t\leq\tau_i^n}}
\def\en{\varepsilon_n}

\def\co{C_0}
\def\supt{\sup_{0\leq t\leq T}}
\def\supn{\sup_{n\geq0}}
\def\pa{\mu}
\def\cvas{\stackrel{a.s.}{\longrightarrow}}
\def\st{\theta_k}
\def\lmin{\lambda_{\mathrm{min}}}
\def\lmax{\lambda_{\mathrm{max}}}
\def\tr{\operatorname{Tr}}
\def\appi{\chi_\mu}

\newcommand\Abs{{(\textbf{A}$_{\sigma}$)}}
\newcommand\Abse{{(\textbf{A}$^{\mathrm{Ellip.}}_{\sigma}$)}}
\newcommand\Au{{(\textbf{A}$_u$)}}
\newcommand\AN{{(\textbf{A}$_N$)}}
\newcommand\AS{{(\textbf{A}$_S$)}}
\newcommand\AY{{(\textbf{A}$_Y$)}}

\newcommand\Diag{\operatorname{Diag}}

\newcommand\1{\mathbf{1}}
\newcommand\dX{\,\mathrm{d} X}
\newcommand\dz{\,\mathrm{d} z}
\newcommand\dD{\partial\D}

\newcommand\cF{{\mathcal F}}
\newcommand\cM{{\mathcal M}}
\newcommand\bN{{\mathbb N}}
\newcommand\cN{{\mathcal N}}
\newcommand\cT{{\mathcal T}}
\newcommand\dr{\,\mathrm{d}r}
\newcommand\dt{\,\mathrm{d}t}
\newcommand\dB{\,\mathrm{d}B}
\newcommand\dM{\,\mathrm{d}M}
\newcommand\dS{\,\mathrm{d}S}
\newcommand\di{ }
\newcommand\ths{\theta_\sigma}
\newcommand\dbi{D_{B,i}}
\newcommand\dsi{D_{S,i}}
\newcommand\li{\Lambda_{\tau^n_{i-1}}}
\newcommand\si{\sigma_{\tau^n_{i-1}}}
\newcommand\rhos{\rho_N}
\makeatother

\begin{document}
\begin{frontmatter}

\title{Almost sure optimal hedging strategy\thanksref{T1}}
\runtitle{Almost sure optimal hedging strategy}

\begin{aug}
\author[a]{\fnms{Emmanuel} \snm{Gobet}\corref{}\ead[label=e1]{emmanuel.gobet@polytechnique.edu}}
\and
\author[b]{\fnms{Nicolas} \snm{Landon}\ead[label=e2]{landon.nico@gmail.com}}
\thankstext{T1}{This research is part of the Chair \textit{Financial Risks}
of the \textit{Risk Foundation}, the Chair \textit{Derivatives of the Future}
and the Chair \textit{Finance and Sustainable Development.}}
\runauthor{E. Gobet and N. Landon}
\affiliation{Ecole Polytechnique and CNRS, and Ecole Polytechnique\break and
CNRS and
GDF SUEZ}
\address[a]{Centre de Mathematiques Appliqu\'ees\\
Ecole Polytechnique and CNRS\\
91128 Palaiseau\\
France\\
\printead{e1}}
\address[b]{GDF SUEZ\\
Direction de la Strat\'{e}gie \\
\quad et du D\'{e}veloppement Durable\\
22, rue du Docteur Lancereaux\\
75008 Paris\\
France\\
and\\
Centre de Mathematiques Appliqu\'ees\\
Ecole Polytechnique and CNRS\\
91128 Palaiseau\\
France\\
\printead{e2}}
\end{aug}

\received{\smonth{1} \syear{2012}}
\revised{\smonth{4} \syear{2013}}

%
\begin{abstract}
In this work, we study the optimal discretization error of stochastic
integrals, in the context of the hedging error in a multidimensional
It\^{o} model when the discrete rebalancing dates are stopping times.
We investigate
the convergence, in an \textit{almost sure} sense, of the renormalized
quadratic variation of the hedging error, for which we exhibit an
asymptotic lower bound for a large class of stopping time strategies.
Moreover, we make explicit a strategy which asymptotically attains this
lower bound \emphh{a.s.} Remarkably, the results hold under great
generality on
the payoff and the model. Our analysis relies on new results enabling
us to control \emphh{a.s.} processes, stochastic integrals and related
increments.
\end{abstract}

%

\end{frontmatter}

\section{Introduction}

\textit{The problem.} We aim at finding a finite sequence of optimal
stopping times $\cT^n=\{\tau^n_0=0 < \tau^n_1 < \cdots < \tau^n_i <
\cdots < \tau^n
_{{N^n_T}
}=T\}$ which minimizes the quadratic variation of the discretization
error of the stochastic integral
\[
Z_s^n=\int_0^s
D_x u(t,S_t) \cdot\dS_t - \sum
_{\tau^n_{i-1}\leq s} D_x u\bigl({\tau^n_{i-1}},S_{\tau^n_{i-1}}
\bigr) \cdot(S_{\tau^n_i\land s }-S_{\tau^n_{i-1}}),
\]
which interpretation is the hedging error \cite{bert:koga:lo:00} of the
discrete Delta-hedging strategy of a European option with underlying
asset $S$ (multidimensional It\^{o} process), maturity $T>0$, price
function $u$ (for the ease of presentation, here $u$ depends only on~$S$)
and payoff $g(S_T)$.
The times $(\tau^n_i)_{1 \leq i \leq{N^n_T}}$ read as rebalancing
dates (or
trading dates), and their number ${N^n_T}$
is a random variable which is finite \emphh{a.s.} The exponent $n$
refers to a
control parameter introduced later on; see Section~\ref{sec2}.
The \emphh{a.s.} minimization of $Z^n_T$ is hopeless since after a suitable
renormalization, it is known that it weakly converges to a mixture of
Gaussian random variables (see \cite{bert:koga:lo:00,GT01,HM05,GT09}
when trading dates are deterministic and under some mild assumptions on
the model and payoff; see \cite{fuka:11} for stopping times under
stronger assumptions). Hence it is more appropriate to investigate the
\emphh{a.s.} minimization of the quadratic variation $\langle
Z^n\rangle_T$
which, owing to the Lenglart inequality (resp., the
Burkholder--Davis--Gundy inequality), allows the control of the
distribution (resp., the $L_p$-moments, $p>0$) of $\sup_{t\leq T
}|Z^n_t|$ under martingale measure. To avoid trivial lower bounds by
letting ${N^n_T}\rightarrow+\infty$, we reformulate our problem into the
\emphh{a.s.} minimization of the product
%
\begin{equation}
\label{eq:critere}{N^n_T} \bigl\langle Z^n \bigr
\rangle_T.
\end{equation}
As emphasized in \cite{F11}, the resolution of this optimization
problem allows the asymptotic minimization of more general costs of the
form $C({N^n_T},\langle Z^n\rangle_T)$, where the function $C\dvtx\R
^2\mapsto\R
$ is increasing in both variables.
Our Theorem \ref{th1} states that the renormalized error \eqref
{eq:critere} has \emphh{a.s.} an asymptotic lower bound over the class of
admissible strategies which consist (roughly speaking\setcounter{footnote}{1}\footnote{A
precise definition is given in Section~\ref{sec2}.}) of deterministic
times and of hitting times of random ellipsoids of the form
%
\begin{equation}
\label{eq:ta:ellip}\hspace*{6pt}\qquad\tau^n_0:=0, \qquad\tau^n_{i}:=
\inf \bigl\{t\geq\tau^n _{i-1}\dvtx(S_t-S_{\tau^n_{i-1}})
\cdot H_{\tau^n_{i-1}} (S_t-S_{\tau^n_{i-1}})= 1 \bigr\} \land T,
\end{equation}
where $(H_t)_{0\leq t \leq T}$ is a measurable adapted
positive-definite symmetric matrix process. It includes the Karandikar
scheme \cite{K95} for discretization of stochastic integrals. In
addition, in Theorems \ref{th2} and \ref{th2:bis} we show the existence
of a strategy of the hitting time form attaining the \emphh{a.s.} lower bound.
The derivation of a central limit-type theorem for $Z^n$ is left to
further research {(see \cite{land:13})}, in particular because the
verification of the criteria in \cite{fuka:11} is difficult to handle
in our general setting.

\textit{Literature background}. Our work extends the existing
literature on discretization errors for stochastic integrals with
deterministic time mesh, mainly considered with financial applications.
Many works deal with hedging rebalancing at regular intervals of length
$\Delta t_i=T / n$. In \cite{Zh99} and \cite{bert:koga:lo:00}, the
authors show that $\E[\langle Z^n \rangle_T]$ converges to 0 at rate
$n$ for payoffs smooth enough [this convergence rate originates to
consider the product \eqref{eq:critere} as a minimization criterion].
However, in \cite{GT01} it is proved that the irregularity of the
payoff may deteriorate the convergence rate: it becomes $n^{1/2}$ for
digital call option. This phenomenon has been intensely analyzed by
Geiss and his co-authors using the concept of fractional smoothness
(see \cite{geis:geis:04,gobe:makh:10,GG11,geis:geis:gobe:11} and
references therein): by the choice of rebalancing dates suitably
concentrated at maturity, we recover the rate $n$.

The first attempt to find optimal strategies with nondeterministic
times goes back to \cite{MP99}: the authors allow a fixed number $n$ of
random rebalancing dates, which actually solve an optimal
multiple-stopping problem. Numerical methods are required to compute
the solution. In \cite{F11}, Fukasawa performs an asymptotic analysis
for minimizing the product $\E({N^n_T})\E(\langle Z^n\rangle_T)$ (an
extension to jump processes has been recently done in \cite
{rose:tank:11}). Under regularity and integrability assumptions (and
for a convex payoff on a single asset), Fukasawa derives an asymptotic
lower bound and provides an optimal strategy. His contribution is the
closest to our current work. But there are major differences:
\begin{longlist}[(1)]
\item[(1)] We focus on \emphh{a.s.} results, which is probably more meaningful
for hedging issues. We are not aware of similar works in this direction.
\item[(2)] We allow a quite general model for the asset. It can be a
multidimensional diffusion process (local volatility model); see the
discussion in Section~\ref{subsection:Au}. As a comparison, in
\cite
{F11} the analysis is carried out for a one-dimensional model (mainly
Black--Scholes model).
\item[(3)] We also allow a great generality on the payoff. In
particular, the payoff can be discontinuous, and the option can be
exotic (Asian, lookback, \ldots) (see Section~\ref{subsection:Au} for
examples): for mathematical reasons, this is a major difference in
comparison with \cite{F11}. Indeed, in the latter reference, the payoff
convexity is needed to ensure the positivity of the option Gamma
(second derivative of price), which is a crucial property in the
analysis. Also, for discontinuous payoff the $L_p$ integrability of the
sensitivities (Greeks) up to maturity may be not satisfied (see \cite
{gobe:makh:11}); thus, some quantities in the analysis (e.g., the
integral of the second moment of the Gamma of digital call option) may
become infinite. In our setting, we circumvent these issues by only
requiring the sensitivities to be finite \emphh{a.s.} up to maturity: actually,
this property is systematically satisfied by payoffs for which the
discontinuity set has a zero-measure (see Section~\ref{subsection:Au}), which includes all the usual situations to our knowledge.
\end{longlist}
To achieve such a level of generality and an \emphh{a.s.} analysis, we design
efficient tools to analyze the \emphh{a.s.} control and \emphh{a.s.}
convergence of local
martingales, of their increments and so forth. All these results
represent another important theoretical contribution of this work.
Other applications of these techniques are in preparation. At last,
although the distribution of hitting time of random ellipsoid of the
form \eqref{eq:ta:ellip} is not explicit, quite surprisingly we obtain
tight estimates on the maximal increments of $\sup_{i\leq
{N^n_T}}(\tau^n
_{i}-\tau^n_{i-1})$, which may have applications in other areas (like
stochastic simulation).

\textit{Outline of the paper}.
In the following, we present some notation and assumptions
that will be used throughout the paper. Section~\ref{sec2} is aimed at
defining our class of stopping time strategies and deriving some
general theoretical properties in this class. For that, we establish
new key results about \emphh{a.s.} convergence, which fit well our framework.
All these results are not specifically related to financial
applications. The main results about hedging error are stated and
proved in Section~\ref{sec3}. Numerical experiments are presented in
Section~\ref{sec4}, with a practical description of the algorithm to
build the optimal sequence of stopping times (actually hitting times)
and a numerical illustration regarding the exchange binary option (in
dimension~2).

\textit{Notation used throughout the paper.}
\begin{itemize}
\item We denote by $x\cdot y$ the scalar product between two vectors
$x$ and $y$, and by $|x|=(x\cdot x)^{1/2} $ the Euclidean norm of $x$;
the induced norm of a $m\times d$-matrix $A$ is denoted by $|A|:=\sup_{x\in\R^d\dvtx|x|=1} |A x|$.
\item$A^*$ stands for the transposition of the matrix $A$; $I_d$
stands for the identity matrix of size $d$; the trace of a square
matrix $A$ is denoted by $\tr(A)$.
\item$\mathcal{S}^d(\R)$, $\mathcal{S}_+^d(\R)$ and $\mathcal
{S}_{++}^d(\R)$ are respectively the set of symmetric, symmetric
nonnegative-definite and symmetric positive-definite $d\times
d$-matrices with coefficients in $\R$: $A\in\mathcal{S}_+^d(\R)$
[resp., $\mathcal{S}_{++}^d(\R)$] if and only if $x\cdot Ax \geq0$
(resp., $>0$) for any $x\in\R^d\setminus\{0\}$.
\item For $A\in\mathcal{S}^d(\R)$, $\Lambda(A):=(\lambda
_1(A),\ldots,\lambda_d(A))$ stands for its spectrum (its $\R$-valued
eigenvalues), and we set $\lmin(A):=\min_{1\leq i\leq d} \lambda_i(A)$.
\item For the partial derivatives of a function $f\dvtx(t,x,y)\mapsto
f(t,x,y)$, we write $D_t f(t,x,y)=\frac{\partial f}{\partial
t}(t,x,y)$, $D_{x_i} f(t,x,y)=\frac{\partial f}{\partial x_i}(t,x,y)$,
$D_{x_i x_j}^2 f(t,x,y)=\break  \frac{\partial^2 f}{\partial x_i\,\partial
x_j}(t,x,y)$, $D^2_{x_iy_j} f(t,x,y)= \frac{\partial^2 f}{\partial
x_i\,
\partial y_j}(t,x,y)$ and so forth.
\item When convenient, we adopt the short notation $f_t$ in place of
$f(t,S_t,Y_t)$ where $f$ is a given function and $(S_t,Y_t)_{0\leq t
\leq T}$ is a continuous time process (introduced below).
\item For a $\R^d$-valued continuous semimartingale $M$, $\langle
M\rangle_t$ stands for the matrix of cross-variations $(\langle
M^i,M^j\rangle_t)_{1\leq i,j\leq d}$.
\item
The constants of the multidimensional version of the
Burkholder--Davis--Gundy inequalities \cite{KS91}, page 166, are defined
as follows: for any $p>0$ there exists $c_p>1$ such that for any vector
$M=(M^1,\ldots,M^d)$ of continuous local martingales with $M_0=0$ and
any stopping time $\theta$, we have
%
\begin{equation}
\label{eq:bdg} c^{-1}_p \E\Biggl|\sum
_{j=1}^d \bigl\langle M^j \bigr
\rangle_{\theta} \Biggr|^{p}\leq\E \Bigl(\sup_{t\leq\theta}
\llvert M_{t}\rrvert^{2p} \Bigr)\leq c_p\E\Biggl|
\sum_{j=1}^d \bigl\langle M^j
\bigr\rangle_{\theta} \Biggr|^{ p}.
\end{equation}
\item For a given sequence of stopping times $\cT^n$, the last time
before $t\leq T$ is defined by $\varphi(t)=\max\{\tau^n_j;\tau
^n_j\leq t\}$:
although dependent on $n$, we omit to indicate this dependency to
alleviate notation.
Furthermore, for a process $(f_t)_{0\leq t \leq T}$, we write $\Delta
f_t:=f_t-f_{\varphi(t-)}$ (omitting again the index $n$ for
simplicity); in particular, we have
$\Delta f_{\tau^n_i}=f_{\tau^n_i}-f_{\tau^n_{i-1}}$. Besides we set
$\Delta
_t=t-\varphi(t-)$ and $\Delta\tau^n_i:=\tau^n_i-\tau^n_{i-1}$.
\item We shortly write $X^n\cvas$ if the random variables
$(X^n)_{n\geq
0}$ converge almost surely as $n\rightarrow\infty$. We write
$X^n\cvas
X^{\infty}$ to additionally indicate that the almost sure limit is
equal to $X^{\infty}$. We shall say that the sequence $(X^n)_{n\geq0}$
is bounded if $\sup_{n\geq0} |X^n|<+\infty$, \emphh{a.s.}
\item$\co$ is a \emphh{a.s.} finite nonnegative random variable,
which may
change from line to line.
\end{itemize}

\textit{Model}.
Let $T>0$ be a given terminal time (maturity), and let $(\Omega,\F,\break (\F
_t)_{0 \leq t \leq T}, \P)$ be a filtered probability space, supporting
a $d$-dimensional Brownian motion $B=(B^i)_{1\leq i\leq d}$ defined on
$[0,T]$, where $(\F_t)_{0\leq t \leq T}$ is the $\P$-augmented natural
filtration of $B$ and $\F=\F_T$. This stochastic basis serves as a
modeling of the evolution of $d$ tradable risky assets without
dividends, which price processes are denoted by $S=(S^i)_{1\leq i\leq
d}$. Their dynamics are given by an It\^o
continuous semimartingale which solves
%
\begin{equation}
\label{eq:S} S_t= S_0 + \int_0^t
b_s \,\mathrm{ d}s+ \int_0^t \sigma
_s \dB_s
\end{equation}
with measurable and adapted coefficients $b$ and $\sigma$.
This is the usual framework of complete market; see \cite{MR05}.
Assumptions on $\sigma$ are given below. Furthermore, for the sake of
simplicity we directly assume that the return of the money market
account $(r_t)_t$ is zero and that $b\equiv0$. This simplification is
not really a restriction (see \cite{MR05} for details): indeed, first
we can still re-express prices in the money market account num\'eraire;
second, because we deal with \emphh{a.s.} results, we can consider dynamics
under any equivalent probability measure, and we choose the martingale measure.

From now on, $S$ is a \emph{continuous local martingale}, and $\sigma$
satisfies the following assumption.

\begin{longlist}
\item[\Abs]
\emphh{a.s.} for any $t\in[0,T]$ $\sigma_t$ is nonzero; moreover $\sigma$ satisfies the continuity condition: there exist
a parameter $\ths\in(0,1]$ and a nonnegative \emphh{a.s.} finite random
variable $\co$ such that
\[
|\sigma_t-\sigma_s|\leq\co \bigl(|S_t-S_s|^{\ths}+|t-s|^{
{\ths }/2}
\bigr)\qquad\forall0\leq s,t\leq T \ \mbox{a.s.}
\]
\end{longlist}
The above continuity condition is satisfied if $\sigma_t:=\sigma
(t,S_t)$ for a function $\sigma(\cdot)$ which is $\ths$-H\"older continuous
w.r.t. the parabolic distance. For some of our results, the above
assumption is strengthened into the following:
\begin{longlist}
\item[\Abse]
Assume \Abs\ and that
$\sigma_t$ is elliptic
in the sense
\[
0<\lmin \bigl(\sigma_t\sigma^*_t \bigr)\qquad \forall 0\leq
t\leq T\ \mbox{a.s.}
\]
\end{longlist}
The assumption \Abse\ is undemanding, since we do not suppose any
uniform (in~$\omega$) lower bound.

We consider an exotic option written on $S$ with payoff $g(S_T,Y_T)$
where $Y_T$ is a functional of $(S_t)_{0\leq t\leq T}$. In the
subsequent asymptotic analysis, we assume that
$Y=(Y^i)_{1\leq i\leq d'}$ is a vector of adapted
continuous
nondecreasing processes.
Examples of such an option are given below: this illustrates that the
current setting covers numerous relevant situations beyond the case of simple
vanilla options [with payoff of form $g(S_T)$].
%
%
\begin{example}
(1) Asian options: $Y^j_t:=\int_0^t S^j_s \,\mathrm{ d}s$ and
$g(x,y):=\break  (\sum_{1\leq j\leq d}\pi_j y^j-K)_+$, for some weights $\pi
_j$ and
a given $K\in\R$.
\begin{longlist}
\item[(2)] Lookback options: $Y^j_t:=\max_{0\leq s\leq t} S^j_s$ and
$g(x,y):=\sum_{1\leq j\leq d}(\pi_j y^j - \pi_j' x^j)$.
\end{longlist}
\end{example}
Furthermore, we assume that the price at time $t$ of such an option is
given by $u(t,S_t,Y_t)$ where $u$ is a $\mathcal{C}^{1,3,1}
([0,T[\,\times\R^d\times\R^{d'} )$ function verifying
%
\begin{eqnarray}
\label{eq:prix:int:sto} u(T,S_T,Y_T)&=&g(S_T,Y_T)
 \quad \mbox{and}
 \nonumber
 \\[-8pt]
 \\[-8pt]
 \nonumber
   u(t,S_t,Y_t) &=& u(0,S_0,Y_0)
+ \int_0^t D_xu
(s,S_s,Y_s) \cdot\dS_s
\end{eqnarray}
for any $t\in[0,T]$. The above set of conditions is related to
probabilistic and analytical properties. First, although not strictly
equivalent, it essentially means that the pair $(S,Y)$ forms a Markov
process and this originates why the randomness of the fair price $\E
(g(S_T,Y_T)|\cF_t)$ at time $t$ only comes from $(S_t,Y_t)$. Observe
that this Markovian assumption about $(S,Y)$ is satisfied in the above
examples. Second, the regularity of the price function $u$ is usually
obtained by applying PDE results thanks to Feynman--Kac
representations: it is known that the expected regularity can be
achieved under different assumptions on the smoothness of the
coefficients of $S$ and $Y$, of the payoff $g$, combined with some
appropriate nondegeneracy conditions on $(S,Y)$. The pictures are
multiple, and it is not our current aim to list all the known related
results; we refer to \cite{wilm:dewy:howi:94} for various Feynman--Kac
representations related to exotic options, and to \cite{pasc:11} for
regularity results and references therein. See Section~\ref{subsection:Au} for extra regularity results.
Besides, we assume
\begin{longlist}
\item[\Au] Let $\A\in\D:= \{D^2_{x_j x_k}, D^3_{x_j x_k x_l},
D^2_{t x_j}, D^2_{x_j y_m} \dvtx1\leq j,k,l \leq d,\break 1\leq m\leq d' \}$,
\[
\P \Bigl(\lim_{\delta\rightarrow0}\sup_{0\leq t< T}\sup
_{|x-S_t|\leq
\delta, |y-Y_t|\leq\delta} \bigl|\A u(t,x,y) \bigr|<+\infty \Bigr)=1.
\]
\end{longlist}
Observe that the above assumption is really weak: this is a pathwise
result, and we do not require any $L_p$-integrability of the
derivatives of $u$. In Section~\ref{subsection:Au}, we provide an
extended list of payoffs (continuous or not) of options (vanilla,
Asian, lookback) in log-normal or local volatility models, for which
\Au
\ holds. Even for the simple option payoff $g(S_T)$ in the simple
log-normal model, we have not been able to exhibit a payoff function
$g$ for which \Au\ is not satisfied.

\section{Class ${\mathcal T}^{\mathrm{adm.}}$ of strategies and
convergence results}
\label{sec2}
In this section, we define the class of strategies under consideration,
and establish some preliminary almost sure convergence results in
connection with this class.

A strategy is a finite sequence of increasing stopping times $\{\tau
_0=0 < \tau_1 < \cdots < \tau_i < \cdots < \tau_{N_T}=T\}$ (with
$N_T<+\infty$ \emphh{a.s.}) which stand for the rebalancing dates.
Furthermore,
the number of risky assets held on each interval $[\tau_i, \tau_{i+1})$
follows the usual Delta-neutral rule $D_xu(\tau_i,S_{\tau_i},Y_{\tau_i})$.

\subsection{Assumptions}
Now to derive \emph{asymptotically} optimal results, we consider a
sequence of strategies indexed by the integers $n=0,1,\ldots,$ that is, writing
\[
\cT^n:= \bigl\{\tau^n_0=0 <
\tau^n_1 < \cdots < \tau^n_i < \cdots
{<} \tau^n_{{N^n_T}}\bigr\} \qquad\mbox{for }n=0,1,\ldots,
\]
and we define an appropriate \emph{asymptotic framework}, as the
convergence parameter~$n$ goes to infinity.
Let $(\en)_{n\geq0}$
be a sequence of positive deterministic real numbers converging to 0 as
$n\rightarrow\infty$; assume that it is a square-summable sequence
%
\begin{equation}
\label{eq:sum:en} \sum_{n\geq0}\en^{2} < +
\infty.
\end{equation}
On the one hand, the parameter $\en^{-2\rhos}$ (for some $\rhos\geq1$)
upper bounds (up to a constant) the number of rebalancing dates of the
strategy $\cT^n$, that is:
\begin{longlist}
\item[\AN]
The following nonnegative random variable is \emphh{a.s.} finite:
\[
\supn \bigl(\en^{2\rhos}{N^n_T} \bigr) < +\infty
\]
for a parameter $\rhos$ satisfying $1\leq\rhos< (1+\frac{\ths
}{2})\land\frac{4}3$.
\end{longlist}
On the other hand, the parameter~$\en$ controls the size of variations
of $S$ between two stopping times in $\cT^n$.
\begin{longlist}
\item[\AS]
The following nonnegative random variable is \emphh{a.s.} finite:
\[
\supn \Bigl(\en^{-2}\supi\sup_{t\in(\tau^n_{i-1},\tau^n_i] }|S_t-
S_{\tau^n
_{i-1}}|^2 \Bigr)<+\infty.
\]
\end{longlist}
Observe that assumptions \AN\ and \AS\  play complementary (and not
equivalent) roles.
We are now ready to define the class of sequence of strategies in which
we are seeking the optimal element.
%
%
\begin{definition}
A sequence of strategies $\cT:=\{\cT^n: n\geq0\}$ is \textit{admissible}
if it fulfills the hypotheses \AN\ and \AS. The set of admissible
sequences $\cT$ is denoted by ${\mathcal T}^{\mathrm{adm.}}$.
\end{definition}
The above definition depends on the sequence $(\varepsilon_n)_{n\geq
0}$, which is fixed from now on.\vadjust{\goodbreak}

%
\begin{remark}
\begin{itemize}
\item The larger $\rhos$, the wider the class of strategies under
consideration. The choice $\rhos=1$ is allowed, but seemingly it rules
out deterministic strategies; see the next remark.
\item If $\rhos>1$, a strategy $\cT^n$ consisting of
${N^n_T}=1+\lfloor
\en
^{-2\rhos}\rfloor$ deterministic times with mesh size $\sup_{1\leq
i\leq
{N^n_T}}\Delta\tau^n_i\leq C \en^{2\rhos}$ (this includes the cases of
uniform and some nonuniform time grids) forms an admissible sequence
of strategies, thanks to the $\frac{ 1 }{2 }^-$-H\"older property of
the Dambis--Dubins--Schwarz Brownian motion of $S^j$ ($1\leq j \leq d$)
(under the additional assumption that $\sigma$ is uniformly bounded to
safely maintain the time-changes into a fixed compact interval).
\item Our setting allows us to consider stopping times satisfying the
\emph{strong predictability condition} (i.e., $\tau^n_i$ is $\cF
_{\tau^n
_{i-1}}$-measurable); see \cite{jaco:prot:12}, Chapter~14.
\item
We show in Proposition \ref{prop4} that the strategy $\cT^n$ of
successive hitting times of ellipsoid of size $\varepsilon_n$ forms a
sequence in ${\mathcal T}^{\mathrm{adm.}}$.
\item In Sections~\ref{subsection:tau}--\ref{subsection:vquad}, we
investigate properties of admissible sequences of strategies. Among
others, we show that the mesh size of $\cT^n$ shrinks \emphh{a.s.} to
0, and we
establish tight \emphh{a.s.} upper bounds (see Corollary \ref{cor2}):
namely for
any $\rho\in(0,2]$, there is a \emphh{a.s.} finite random variable
$C_\rho$ such that
$\sup_{1\leq i\leq{N^n_T}} \Delta\tau^n_i\leq C_\rho\en^{2-\rho
}$ for any
$n\geq0$.
\end{itemize}
\end{remark}

We require an extra technical condition on the nondecreasing process
$Y$ which is fulfilled in practical cases for an admissible sequence of
strategies.
\begin{longlist}
\item[\AY] The following nonnegative random variable is \emphh{a.s.} finite:
for some $\rho_Y>4(\rhos-1)$
\[
\supn \Bigl(\en^{-\rho_Y}\supi|\Delta Y_{\tau^n_i}| \Bigr)<+\infty.
\]
\end{longlist}

%
\begin{example} Let $\cT:=\{\cT^n\dvtx n\geq0\}$ satisfy \AS--\AN.
\begin{longlist}[(1)]
\item[(1)] Asian options: applying Corollary \ref{cor2}
[item (ii)] with {$\rho=\frac{2}3$ and taking $\rho_Y=\frac{4}3> 4(\rhos-1)$
(since $\rhos<\frac{ 4 }{3 }$)} gives
\[
\supn \Bigl({\en^{-\rho_Y}}\supi|\Delta Y_{\tau^n_i}| \Bigr)\leq \sup
_{0\leq t
\leq T}|S_t| \supn \Bigl({\en^{\rho-2}}\supi
\Delta\tau^n_i \Bigr)<+\infty\qquad \mbox{a.s.}
\]
\item[(2)] Lookback
options: clearly, we have
\[
\supn \Bigl(\en^{-1}\supi|\Delta Y_{\tau^n_i}| \Bigr)\leq\supn
\Bigl(\en^{-1}\supt|\Delta S_t| \Bigr)<+\infty\qquad\mbox{a.s.};
\]
{thus \AY\ is satisfied with $\rho_Y=1$ provided that $\rhos<5/4$.}
\end{longlist}
\end{example}

\subsection{Fundamental lemmas about almost sure convergence}
This subsection is devoted to the main ingredient (Lemmas \ref{lem1}
and \ref{lem1:bis}) about almost sure convergence, which is involved in
the subsequent asymptotic analysis.

We first recall some usual approaches to establish that a sequence
$(U^n_T)_{n\geq0}$ converges to 0 in probability or almost surely, as
$n\rightarrow\infty$: it serves as a preparation for the comparative
discussion we will have regarding our almost sure convergence results.
\begin{itemize}
\item\emph{Convergence in probability.}
It can be handled, for instance, by using the Markov inequality and
showing that the $L_p$-moment (for some $p>0$) of $U^n_T$ converges to
0: for $p=1$ and $\delta>0$, it writes $\P(|U^n_T|\geq\delta)\leq
\frac
{ \E|U^n_T| }{ \delta}\rightarrow_{n\rightarrow\infty}0$. Observe that
this approach requires a bit of integrability of the random variable
$U^n_T$.

To achieve the uniform convergence in probability of $(U^n_{t})_{0\leq
t \leq T}$ to 0, Lenglart \cite{leng:77} introduced an extra condition:
the relation of domination. Namely, assume that $(U^n_t)_{0\leq t \leq
T}$ is a nonnegative continuous adapted process and that it is
dominated by a nondecreasing continuous adapted process
$(V^n_t)_{0\leq t \leq T}$ (with $V^n_0=0$) in the sense $\E
(U^n_\theta
)\leq\E(V^n_\theta)$ for any stopping time $ \theta\in[0,T]$. Then,
for any $c_1, c_2>0$, we have
\[
\P \Bigl(\sup_{t\leq T} U^n_t\geq
c_1 \Bigr)\leq\frac{ 1 }{c_1 }\E \bigl(V^n_T
\land c_2 \bigr)+\P \bigl(V^n_T\geq
c_2 \bigr).
\]
A standard application consists in taking $U^n$ as the square of a
continuous local martingales $M^n$; then, the convergence in
probability of $\langle M^n,M^n\rangle_T$ to 0 implies the uniform
convergence in probability of $(M^n_{t})_{0\leq t \leq T}$ to 0. The
converse is also true, the relation of domination deriving from BDG
inequalities. This kind of result leads to useful tools for
establishing the convergence in probability of triangular arrays of
random variables: for instance, see \cite{geno:jaco:93}, Lemma 9, in the
context of parametric estimation of stochastic processes.
\item\emph{Almost sure convergence.}
We may use a Borel--Cantelli type argument, assuming that $\sum_{n\geq
0}\E
|U^n_T|<+\infty$. Fubini--Tonelli's theorem yields that the series
$\sum_{n\geq0}|U^n_T|$ converges a.s., and in particular $U^n_T\cvas
0$. Here
again, the integrability of $U^n_T$ is required.

Bichteler and Karandikar
leveraged this type of series argument to establish the \emphh{a.s.} convergence
of stochastic integrals under various assumptions, with in view either
approximation issues or pathwise stochastic integration; see \cite
{B81,K89,K95,K06} and references therein.
\end{itemize}
Our result below (Lemma \ref{lem1}) is inspired by the above
references, but its conditions of applicability are less stringent, and
it allows more flexibility in our framework.
We assume a relation of domination, but:
\begin{longlist}[(1)]
\item[(1)] not for all stopping times (as in Lenglart domination);
\item[(2)] the processes $(U^n_t)_{0\leq t\leq T}$ are not assumed to
be continuous
[nor\break  $(\sum_{n\geq0}U^n_t)_{0\leq t\leq T}$];
\item[(3)] the dominating process $V^n$ is not assumed to be nondecreasing.
\end{longlist}
Thus, our assumptions are less demanding, but on the other hand, we do
not obtain any uniform convergence result.
Moreover, we emphasize that we do not assume any integrability on
$U^n_T$. This is crucial, because the typical applications of Lemma~\ref
{lem1} are related to $U^n_T$ defined as a (possibly stochastic)
integral of the derivatives of $u$ evaluated along the path
$(S_t,Y_t)_{0\leq t \leq T}$: since usual payoff functions are
irregular, it is known that the $L_p$-moments of related derivatives
blow up as time goes to maturity, and it is hopeless to obtain the
required integrability on $U^n_T$ assuming only \Au.

We are now ready for the statement of our \emphh{a.s.} convergence result.
%
%
\begin{lemma}
\label{lem1}
Let $\cM_0^+$ be the set of nonnegative measurable processes vanishing
at $t=0$.
Let $(U^n)_{n\geq0}$ and $(V^n)_{n\geq0}$ be
two sequences of processes in $\cM_0^+$.
Assume that:
\begin{longlist}[(iii)]
\item[(i)] the series $\di\sum_{n\geq0}V^n_t$ converges for all
$t\in[0,T]$,
almost surely;
\item[(ii)]
the above limit is upper bounded by a process $\bar V\in\cM_0^+$ and
that $\bar V$ is continuous a.s.;
\item[(iii)] there is a constant $c\geq0 $ such that,
for every $n\in\bN$, $k\in\bN$ and $t\in[0,T]$, we have
\[
\E \bigl[U^n_{t\wedge\st} \bigr]\leq c \E \bigl[V^n_{t \wedge\st}
\bigr]
\]
with the random time $\st:=\inf\{s\in[0,T]\dvtx\bar V_s \geq k\}
$.\footnote{With the usual convention $\inf\varnothing=+\infty$.}
\end{longlist}
Then for any $t\in[0,T]$, the series $\di\sum_{n\geq0}U^n_t$ converges
almost surely. As a consequence, $U^n_t\cvas0$.
\end{lemma}

\begin{pf}
First, observe that $(\st)_{k\geq0}$ defines well random times since
$\bar V$ is continuous.

Denote by $\cN_V$ the $\P$-negligible set on
which the series $(\sum_{n\geq0}V^n_t)_{0\leq t \leq T}$ do not
converge, and
on which $\bar V$ and then $(\st)_{k \geq0}$ are not defined; observe
that for $\omega\notin\cN_V$, we have $\bar V_{t\wedge\st}(\omega
)\leq k$ for any $t\in[0,T]$ and $k\in\bN$.
Set $\bar V^p:= \sum_{n=0}^p V^n$: we have $\bar V^p\leq\bar V$ on
$\cN_V^c$; thus, the localization of $\bar V$ entails that of $\bar
V^p$ and we have $\bar V^p_{t \wedge\st}\leq k$ for any $k, p$ and $t$
(on $\cN_V^c$).\vspace*{1pt}

Moreover, for any $n$ and $k$, the relation of domination
writes
%
\begin{equation}
\label{eq:local:proof} \E \Biggl[\sum_{n=0}^p
U_{t \wedge\st}^n \Biggr]\leq c \E \Biggl[\sum
_{n=0}^p V_{t \wedge\st}^n \Biggr]=
c\E \bigl[ \bar V^p_{t \wedge
\st} \bigr]\leq c k.
\end{equation}
From Fatou's lemma, we get $\E[\sum_{n\geq0}U_{t \wedge\st}^n
]<+\infty$:
in particular, for any $k\in\mathbb{N}$, there is a $\P$-negligible set
$\mathcal{N}_{k,t}$,
such that $\sum_{n\geq0}U_{t \wedge\st}^n(\omega)$ converges for all
$\omega
\notin\mathcal{N}_{k,t}$. The set $\cN_t=\bigcup_{k\in\mathbb{N}}
\mathcal{N}_{k,t}\cup\cN_V$ is $\P$-negligible,
and it follows that for $\omega\notin\cN_t$, the series $\sum_{n\geq0}
U^n_{t\wedge\st}(\omega)$ converges for all $k\in\bN$. For $\omega
\notin\cN_t$, we have $\theta_k(\omega)=+\infty$ as soon as
$k>\bar
V_T(\omega)$; thus by taking such $k$, we complete the convergence of
$\sum_{n\geq0}U^n_t$ on $\cN_t^c$.
\end{pf}
Observe that in our argumentation, we do not assume that the
nonnegative random variables $U^n_t$ and $V^n_t$ have a finite
expectation (and in some examples, it is false, especially at $t=T$).
However, note that in \eqref{eq:local:proof} we prove that
$U^n_{t\wedge\st}$ and $V^n_{t\wedge\st}$ have a finite expectation:
in other words, $(\st)_{k\geq0}$ serves as a common localization for
$U^n$ and $V^n$. In addition, Lemma \ref{lem1} is general and thorough
since we do not assume any adaptedness or regularity properties of the
processes $U^n$ and $V^n$. We provide a simpler version that can be
customized for our further applications:
%
%
\begin{lemma}
\label{lem1:bis} Let $\C_0^+$ be the set of nonnegative continuous
adapted processes, vanishing at $t=0$.
Let $(U^n)_{n\geq0}$ and $(V^n)_{n\geq0}$ be two sequences of
processes in $\C_0^+$. Replace the two first items of Lemma \ref
{lem1} by:
\begin{longlist}[(iii$'$)]
\item[(i$'$)] $t\mapsto V^n_t$ is a nondecreasing function on $[0,T]$,
almost surely;
\item[(ii$'$)] the series $\di\sum_{n\geq0}V^n_T$ converges almost surely;
\item[(iii$'$)] there is a constant $c\geq0 $ such that,
for every $n\in\bN$, $k\in\bN$ and $t\in[0,T]$, we have
%
\begin{equation}
\E \bigl[U^n_{t\wedge\st} \bigr]\leq c \E \bigl[V^n_{t \wedge\st}
\bigr] \label{eq:rel:dom}
\end{equation}
with the stopping time $\st:=\inf\{s\in[0,T]\dvtx\bar V_s \geq k\}$
setting $\di\bar V_t= \sum_{n\geq0}V^n_t$.
\end{longlist}
Then, the conclusion of Lemma \ref{lem1} still holds.
\end{lemma}
\begin{pf}
We just have to prove that items (i$'$)${} + {}$(ii$'$) entails items
{(i)${} + {}$(ii)} of Lemma \ref{lem1} for $U^n$ and $V^n$ in $\C
_0^+\subset
\cM_0^+$. Since $V^n$ is nondecreasing, the \emphh{a.s.} convergence
of $\sum_{n\geq0}
V^n_T$ implies that of $\sum_{n\geq0}V^n_t$. Moreover $\sum_{n\geq
0}\sup_{0\leq
t\leq T} V_t^n = \sum_{n\geq0}V_T^n <+\infty$ \emphh{a.s.} Therefore,
\emphh{a.s.} the
series associated with $V^n$ is normally convergent on $[0,T]$ and
$\bar V:=\sum_{n\geq0}V^n \in\C^+_0$: items (i)${} + {}$(ii) are
satisfied.
Observe $\theta_k$ is a \emph{stopping} time since $\bar V$ is
continuous and adapted.
\end{pf}

We apply Lemma \ref{lem1:bis} to derive a simple criterion for the
convergence of continuous local martingales.
%
%
\begin{corollary}\label{cor1} Let $p>0$, and let $\{(M^n_t)_{0\leq t
\leq T}\dvtx n\geq0\}$ be a sequence of scalar continuous local martingales
vanishing at zero. Then
\[
\sum_{n\geq0} \bigl\langle M^n \bigr
\rangle_T^{p/2} \cvas\quad\Longleftrightarrow\quad\sum
_{n\geq0} \supt\bigl|M^n_t\bigr|^p
\cvas.
\]
\end{corollary}
\begin{pf} We first prove the implication $\Rightarrow$. Set $U^n_t:=
\sup_{0\leq s\leq t} |M^n_s|^p$ and $V^n_t:=\langle M^n\rangle
_t^{p/2}$, and let us check the conditions of Lemma \ref{lem1:bis}:
(i$'$) $V^n$ is nondecreasing and (ii$'$) $\sum_{n\geq0}V^n_T$
converges \emphh{a.s.}
The relation of domination \eqref{eq:rel:dom} follows from the BDG
inequalities [see the RHS of \eqref{eq:bdg}]
and we are done.
The implication $\Leftarrow$ is proved similarly,
using the LHS of \eqref{eq:bdg} regarding the BDG inequalities.
\end{pf}

\subsection{Controls of \texorpdfstring{$\Delta\tau^n$}{Delta tau n} and of the martingales increments}
\label{subsection:tau}
Being inspired by the scaling property of Brownian motion, we might
intuitively guess that a sequence of strategy $(\cT^n)_{n\geq0}$
satisfying \AS\ yields stopping times increments of magnitude equal
roughly to $\varepsilon_n^2$. Actually, thorough estimates are
difficult to derive: for instance, the exit times of balls by a
Brownian motion define unbounded random variables.

To address these issues, we take advantage of Lemma \ref{lem1:bis} to
establish estimates on the sequence $(\Delta\tau^n_i:=\tau^n_i-\tau^n
_{i-1})_{1\leq i \leq{N^n_T}}$, which show that we almost recover the
familiar scaling $\en^2$.
%
%
\begin{proposition}
\label{prop1} Assume \textup{\Abs}. Let $\cT$ be a sequence of strategies
satisfying \textup{\AS}\
and let $p\geq0$. Then:\vadjust{\goodbreak}
\begin{longlist}[(ii)]
\item[(i)] The series $\di\sum_{n\geq0}\en^{-(p-2)}\supi(\Delta
\tau^n_{i})^p
\cvas$.
\item[(ii)] Assume moreover that $\cT\in{\mathcal T}^{\mathrm
{adm.}}$: the series $\di
\sum_{n\geq0}
{\en^{-2(p-1)+2\rhos}}\times \sum_{\tau_{i-1}^n< T}(\Delta\tau^n_{i})^p
\cvas$.
\end{longlist}
\end{proposition}
{The proof is postponed to Appendix \ref{appendix:prop1}.}
As a consequence of Proposition \ref{prop1}, the mesh size of $\cT^n$,
that is, $\sup_{1\leq i\leq{N^n_T}} \Delta\tau^n_i$, converges
\emphh
{a.s.} to 0 as
$n\rightarrow\infty$, with some explicit rates of convergence: this is
the statement below.

%
\begin{corollary}
\label{cor2} With the same assumptions and notation as Proposition~\ref
{prop1}, we have the following estimates, for any $\rho>0$:
\begin{longlist}[(ii)]
\item[(i)] Under \textup{\AS},
$\di\sup_{n\geq0} (\en^{\rho-1}\sup_{1\leq i\leq{N^n_T}} \Delta
\tau^n_i
)<+\infty$ \emphh{a.s.}
\item[(ii)] Under \textup{\AS}--\AN, $\di\sup_{n\geq0} (\en^{\rho
-2}\sup_{1\leq i\leq{N^n_T}} \Delta\tau^n_i )<+\infty$ \emphh{a.s.}
\end{longlist}
\end{corollary}

\begin{pf} Item {(i)}. Clearly, from Proposition \ref{prop1}{(i)}, we obtain\break
$\di\sup_{n\geq0} (\en^{-(p-2)}\times  \sup_{1\leq i\leq
{N^n_T}} (\Delta\tau^n_i)^p )<+\infty$ \emphh{a.s.} for
any $p\geq
0$ and the result follows by taking $p=2/\rho$.

Item {(ii)}. We proceed similarly by observing that Proposition
\ref{prop1}{(ii)} gives
\begin{eqnarray*}
\di\sup_{n\geq0} \Bigl({\en^{-2(p-1-\rhos)}}\sup
_{1\leq i\leq{N^n_T}} \bigl(\Delta\tau^n_i
\bigr)^p \Bigr)&\leq&\sup_{n\geq0} \biggl( {
\en^{-2(p-1-\rhos)}}\sum_{\tau_{i-1}^n<
T}\bigl(\Delta
\tau^n_{i}\bigr)^p \biggr)\\
& <&+\infty\qquad
\mbox{a.s.}
\end{eqnarray*}
\upqed\end{pf}\eject

We are now in a position to control the \emphh{a.s.} convergence of some
stochastic integrals appearing in our further optimality analysis. The
following proposition and corollary will play a crucial role in the
estimations of the error terms appearing in the main theorems; see
Section~\ref{sec3}.

%
\begin{proposition}
\label{prop2}
Assume \textup{\Abs}. Let $\cT=(\cT^n)_{n\geq0}$ be a sequence of strategies,
$((M^n_t)_{0\leq t \leq T})_{n\geq0}$ be a sequence of $\R$-valued
continuous local martingales such that $\langle M^n\rangle_t=\int_{0}^t
\alpha^n_r \dr$ for a nonnegative measurable adapted $\alpha^n$
satisfying the following inequality: there exists a nonnegative \emphh{a.s.}
finite random variable $C_\alpha$ and a parameter $\theta\geq0$ such that
\[
0\leq\alpha^n_r\leq C_\alpha \bigl(|\Delta
S_r|^{2\theta}+|\Delta r|^\theta \bigr)\qquad \forall 0
\leq r< T, \forall n\geq0,\ \mbox{a.s.}
\]
Then, the following convergences hold:
\begin{longlist}[(ii)]
\item[(i)] Assume $\cT$ satisfies \textup{\AS}\ and let $p\geq2$
\[
\di\sum_{n\geq0} \biggl(\en^{3-((1+\theta)/2)p}\sum
_{\tau
_{i-1}^n< T}\supit\bigl|\Delta M^n_t\bigr|^p
\biggr)<+\infty\qquad\mbox{a.s.}
\]
\item[(ii)] Assume furthermore that $\cT$ satisfies \textup{\AN}\ (i.e., $\cT
\in
{\mathcal T}^{\mathrm{adm.}}$), and let $p>0$
\[
\di\sum_{n\geq0} \biggl({\en^{2-(1+\theta)p+2\rhos}}\sum
_{\tau
_{i-1}^n< T}\supit\bigl|\Delta M^n_t\bigr|^p
\biggr)<+\infty\qquad\mbox{a.s.}
\]
\end{longlist}
\end{proposition}
{The proof is postponed in Appendix \ref{appendix:prop2}.}
A straightforward consequence of the aforementioned proposition is
given by the following corollary, which proof is left to the reader.
%
%
\begin{corollary}
\label{cor3}
Using the assumptions and notation of Proposition \ref{prop2}, we have
the following estimates, for any $\rho>0$:
\begin{longlist}[(ii)]
\item[(i)] Under \textup{\AS}, $\di\supn(\en^{\rho-{(1+\theta)
}/2}\supi
\supit|\Delta{M^n_t}| )<+\infty$, \emphh{a.s.}
\item[(ii)] Under \textup{\AS--\AN}, $\di\supn(\en^{\rho-(1+\theta
)}\supi
\supit|\Delta{M^n_t}| )< +\infty$, \emphh{a.s.}
\end{longlist}
\end{corollary}
%
%
\begin{remark}
Observe that in the proofs of the Section~\ref{subsection:tau}
{results}, we have not used the knowledge of the upper bound on $\rhos$
[stated in \AN]: it means that all the related results are true for any
admissible sequence of strategies assuming only $\rhos\geq1$.
\end{remark}
%
\subsection{Almost sure convergence of weighted discrete quadratic variation}
\label{subsection:vquad}
%
%
\begin{proposition}
\label{prop3} Assume \textup{\Abs}\ and let $\cT$ be a sequence of strategies
satisfying \textup{\AS}. Let $(H_t)_{0\leq t < T}$ be a continuous adapted
$d\times d$-matrix process such that $\sup_{t\in[0,T)}|H_t|<+\infty$ a.s.,
and let $(M_t)_{0\leq t \leq T}$ be a $\R^d$-valued continuous local
martingale such that $\langle M\rangle_t=\int_{0}^t \alpha_r \dr$ with
$\supt|\alpha_t|<+\infty$ \emphh{a.s.} Then
\[
\sum_{\tau_{i-1}^n< T}\Delta M_{\tau^n_i}^*H_{\tau^n_{i-1}}
\Delta M_{\tau^n_i} \cvas \int_0^T \tr
\bigl(H_t \,\mathrm{ d}\langle M\rangle_t \bigr).
\]
\end{proposition}

\begin{pf}
From It\^o's lemma, $\sum_{\tau_{i-1}^n< T}\Delta M_{\tau
^n_i}^*H_{\tau^n_{i-1}} \Delta
M_{\tau^n
_i}$ is equal to
\begin{eqnarray*}
&&\sum_{k,l=1}^d \sum
_{\tau_{i-1}^n< T}\Delta M^k_{\tau^n_i}H^{k,l}_{\tau^n_{i-1}}
\Delta M^l_{\tau^n_i}\\
&&\qquad=\sum_{k,l=1}^d
\int_0^T H^{k,l}_{\varphi(t)}
\bigl(\Delta M^k_{t}\dM^l_{t}+
\Delta M^l_{t}\dM^k_{t} +\,\mathrm{ d} \bigl\langle M^k,M^l \bigr\rangle_t
\bigr)
\\
&&\qquad= \int_0^T \Delta M^*_{t}
\bigl(H_{\varphi(t)}+H^*_{\varphi(t)} \bigr) \dM_t +\int
_0^T \tr \bigl(H_{\varphi(t)}\,\mathrm{ d}
\langle M\rangle_t \bigr).
\end{eqnarray*}
The second term in the above RHS converges \emphh{a.s.} to $\int_0^T
\tr
(H_{t}\,\mathrm{ d}\langle M\rangle_t)$: indeed, the difference is
bounded by
$C_0\int_0^T |H_t-H_{\varphi(t)}|\dt$, and we conclude by an
application of the dominated convergence theorem, invoking the
continuity and boundedness of $H$ and the convergence to 0 of the mesh
size of $\cT^n$; see Corollary \ref{cor2}.

Thus it remains to show that the stochastic integral w.r.t. $\dM_t$
converges \emphh{a.s.} to 0. Owing to Corollary \ref{cor1}, it is
enough to
study the series of quadratic variations, that is, to show that
$\sum_{n\geq0}[\int_0^T (\Delta M_t^*(H_{\varphi(t)}+H^*_{\varphi
(t)}) \,\mathrm{ d}\langle M \rangle_t (H_{\varphi(t)}+H^*_{\varphi(t)})
\Delta
M_t ) ]^3\cvas$,
and since $\alpha$ and $H$ are \emphh{a.s.} bounded on $[0,T)$, it is
sufficient
to show
%
\begin{equation}
\label{eq:conv:1}\sum_{n\geq0} \biggl[\int
_0^T | \Delta M_t|^2
\dt \biggr]^3\cvas.
\end{equation}
Clearly $ [\int_0^T |\Delta M_t|^2\dt]^3$ is bounded by
\[
d^3T^3\sup_{1\leq j\leq d} \supi\supit\bigl|\Delta M_t^j\bigr|^6\leq C_0 \en^2
\]
owing to Corollary \ref{cor3} [item (i)] for $\theta=0$ and $\rho
=\frac{1}6$.
The convergence \eqref{eq:conv:1} is proved, and we are done.
\end{pf}

\subsection{Verification of the hypothesis on a special family of
hitting times}
One of the more appealing results of the paper is that a very large
family of hitting times fulfills the assumptions \AN\ and \AS\ with a
threshold depending of $\en$.
%
%
\begin{proposition}
\label{prop4} Assume \textup{\Abs}.
Let $(H_t)_{0\leq t {<} T}$ be a continuous adapted
nonnegative-definite $d\times d$-matrix process, such that \emphh{a.s.}
\[
0<\inf_{0\leq t < T}\lmin(H_t)\leq\sup
_{0\leq t < T}\lmax(H_t)<+\infty.
\]
The strategy $\cT^n$ given by
\[
\cases{ \tau^n_0:=0,
\vspace*{2pt}\cr
\displaystyle\tau^n_{i}:=\inf \bigl\{t\geq\tau^n_{i-1}
\dvtx(S_t - S_{\tau^n_{i-1}})^* H_{\tau^n
_{i-1}} (S_t -
S_{\tau^n_{i-1}}) > \en^2 \bigr\}\wedge T, }
\]
defines a sequence of strategies satisfying assumptions \textup{\AN}\
{[with\break
$\supn(\en^{2}{N^n_T}) < +\infty$ \emphh{a.s.}]} and \textup{\AS}, that is
$\{
\cT^n\dvtx n\geq
0\}\in{\mathcal T}^{\mathrm{adm.}}$.
\end{proposition}
The proof is postponed to Appendix \ref{appendix:lem6}. Observe that
the above sequence of strategies is admissible even in the most
constrained case $\rhos=1$. As we shall see later on, the optimal
stopping times are given by the hitting times by the process $S$ of an
ellipsoid (corresponding to the case $H$ symmetric).\vadjust{\goodbreak}

\section{Main results}
\label{sec3}
\subsection{Statements}
We now go back to the hedging issue: at time $s\in[0,T]$, the fair
value of the option is $u(s,S_s)$, and the hedging portfolio with
discrete rebalancing dates $\cT^n$ is $u(0,S_0)+\sum_{\tau
^n_{i-1}\leq s}
D_x u({\tau^n_{i-1}},S_{\tau^n_{i-1}})\cdot(S_{\tau^n_i\land s
}-S_{\tau^n
_{i-1}})$, which yields an hedging error equal to
%
\begin{eqnarray}\label{eq:zs}
Z_s^n&:=&u(s,S_s)-
\biggl(u(0,S_0)+\sum_{\tau^n_{i-1}\leq s}
D_x u\bigl({\tau^n_{i-1}},S_{\tau^n_{i-1}}
\bigr)\cdot(S_{\tau^n_i\land s
}-S_{\tau^n
_{i-1}}) \biggr)
\nonumber
\\[-8pt]
\\[-8pt]
\nonumber
&=& \int_0^s (D_x
u_t- D_x u_{\varphi(t)} )\cdot\dS_t
\end{eqnarray}
using \eqref{eq:prix:int:sto}, where the integrand appears as the
difference of Delta between $\tau^n_{i-1}$ and $t\in\,]\tau
^n_{i-1},\tau^n_i]$
for each $0\leq i\leq{N^n_T}$.

One main result of the paper is a lower bound of the renormalized
quadratic variation of the hedging error $Z^n$: it is partly derived
from a \emph{smart} representation of
%
\begin{equation}
\label{eq:qvZT} \bigl\langle Z^n\bigr\rangle_T=\int
_0^T (D_x u_t-
D_x u_{\varphi(t)} )^* \,\mathrm{ d}\langle S
\rangle_t (D_x u_t- D_x
u_{\varphi(t)} )
\end{equation}
as a sum of squared random variables and an application of the
Cauchy--Schwarz inequality. To derive this suitable representation, we
apply the It\^o formula and identify the bounded variation term; it is
straightforward in dimension one, much more intricate in a
multidimensional setting and this is equivalent to solving the
following matrix equation.

%
\begin{lemma}
\label{lem5}
Let $c\in\mathcal{S}^d(\R)$. Then the equation
%
\begin{equation}
\label{eq5} 2 \tr(x)x+4x^2=c^2
\end{equation}
admits exactly one solution $x(c)\in\mathcal{S}^d_{+}(\R)$. In
addition, $x(c)$ is positive-definite if and only if $c^2$ is
positive-definite. Last, the mapping $c\mapsto x(c)$ is continuous.
\end{lemma}
The proof is given in Section~\ref{subsection:proof:lem5}.
We are now in a position to give an explicit asymptotic lower bound for
${N^n_T}\langle Z^n\rangle_T$: this is the contents of the following theorem.
%
%
\begin{theorem}
\label{th1}
Assume assumptions \textup{\Abs}, \textup{\Au}, \textup{\AS}, \textup{\AN}\ and \textup{\AY}\ are in force. Let
$X$ be the solution of \eqref{eq5} with $c:=\sigma^* D^2_{xx} u
\sigma$. Then
\[
\liminf_{n\rightarrow+\infty}\ {N^n_T}\bigl\langle
Z^n\bigr\rangle_T\geq \biggl( \int_0^T
\tr(X_t )\dt \biggr)^2 \qquad\mbox{a.s.}
\]
\end{theorem}
Let us comment a bit on the above lower bound:
\begin{itemize}
\item First, it is \emphh{a.s.} finite: indeed, $\sup_{t<T}|\sigma_t^*D_{xx}^2
u_t \sigma_t|<+\infty$ a.s., and the continuity of $c\mapsto x(c)$
imply $\sup_{t<T}|X_t|<+\infty$ \emphh{a.s.}
\item Second, observe that \emphh{a.s.}
\begin{eqnarray*}
&&\biggl\{\int_0^T \tr(X_t) \dt=0
\biggr\} = \bigl\{\forall t<T\dvtx\sigma_t^*D_{xx}^2
u_t\sigma_t=0 \bigr\}
\\
&&\hspace*{75pt}\stackrel{\scriptsize{\mathrm{under}\ \textup{\Abse}}} {=} \bigl\{\forall t<T\dvtx D_{xx}^2
u_t=0 \bigr\}
\end{eqnarray*}
using at the first equality that $\tr(x(c))>0\Leftrightarrow x(c)\neq0
\Leftrightarrow c\neq0$. Then we obtain that except in degenerate
situations [where the Gamma matrix $D^2_{xx} u_t$ is zero at any time,
assuming \Abse], the lower bound in Theorem \ref{th1} is nonzero.
\item As a consequence, we immediately obtain a lower bound for the
$L_p$-criterion: indeed, using the Fatou lemma and the Cauchy--Schwarz
inequality, we derive (for any $p>0$)
\begin{eqnarray*}
\biggl[\E \biggl(\int_0^T
\tr(X_t)\dt \biggr)^p \biggr]^2&\leq& \Bigl[\E
\Bigl(\liminf_{n\rightarrow+\infty} \bigl( {N^n_T}
\bigl\langle Z^n\bigr\rangle _T\bigr)^{p/2}
\Bigr) \Bigr]^2
\\
&\leq&\liminf_{n\rightarrow+\infty} \bigl[\E\bigl({N^n_T}
\bigl\langle Z^n\bigr\rangle_T\bigr)^{p/2}
\bigr]^2
\\
&\leq&\liminf_{n\rightarrow+\infty}\E \bigl(\bigl({N^n_T}
\bigr)^p \bigr)\E \bigl( \bigl\langle Z^n\bigr
\rangle_T^p \bigr).
\end{eqnarray*}
For $p=1$ we recover the Fukasawa approach \cite{F11}.
\end{itemize}
The next theorem tells us that along a suitable sequence $\cT^n$ (the
hitting times of some random ellipsoids) the lower bound {of Theorem
\ref{th1}} is reached. Let $\chi(\cdot)$ be a smooth function such that
$\mathbf{1}_{]-\infty,1/2]}\leq\chi(\cdot)\leq\mathbf{1}_{]-\infty,1]}$
and for $\pa>0$, set $\appi(x)=\chi(x/\mu)$.

%
\begin{theorem}
\label{th2}
Assume assumptions \textup{\Abse}, \textup{\Au}, \textup{\AS}, \textup{\AN}\ and \textup{\AY}\ are in force.
Let $\pa>0$, for $t\geq0$ set $\Lambda_t:=(\sigma^{-1}_t)^*X_t
\sigma^{-1}_t$
and $\Lambda^\pa_t:=\Lambda_t + \pa\appi(\lmin(\Lambda
_t))I_d$.

For a given $n\in\bN$, define the strategy $\cT^n_\pa$ by
%
\begin{equation}
\label{eq1} %
\cases{ \tau^n_0:=0,
\vspace*{2pt}\cr
\displaystyle\tau^n_{i}=\inf \bigl\{t\geq\tau^n_{i-1}
\dvtx(S_t - S_{\tau^n_{i-1}})^* \Lambda^\pa_{\tau^n_{i-1}}
(S_t - S_{\tau^n_{i-1}}) > \en^2 \bigr\}\wedge T.}
\end{equation}
Then, the sequence of strategies $\cT_\pa=\{\cT^n_\pa\dvtx n\geq0\}
$ is
admissible, and it is $\pa$-asymptotically optimal in the following sense:
\[
\limsup_{ n\rightarrow+\infty} \biggl|N_T^{n}\bigl\langle
Z^n\bigr\rangle_T- \biggl( \int_0^T
\tr(X_t) \dt \biggr)^2\biggr |\leq C_\mu\mu\int
_0^T \appi \bigl(\lmin(\Lambda_t)
\bigr)\tr \bigl(\sigma_t \sigma_t^* \bigr)\dt,
\]
where the random variable $C_\pa:= \int_0^T (4 \tr(X_t) + 3\pa
\appi
(\lmin(\Lambda_t))\tr(\sigma_t\sigma_t^*) )\dt$ is \emphh
{a.s.} finite
(locally uniformly w.r.t. $\mu\geq0$).

In particular, on the event $\{\forall t\in[0,T]\dvtx\lmin(\Lambda
_{t})\geq
\pa\}$, $N_T^{n}\langle Z^n\rangle_T$ converges \emphh{a.s.} to $ (
\int_0^T
\tr(X_t)\dt)^2$.
\end{theorem}
Observe that we require the ellipticity condition to hold.
The proof is given in Section~\ref{subsection:proof:th2}.

We can strengthen the above theorem by allowing $\pa=0$ under stronger
assumptions.
%
%
\begin{theorem}
\label{th2:bis}
Assume the assumptions of Theorem \ref{th2} and additionally that
%
\begin{equation}
\label{eq:th3:ass} \P \Bigl(\inf_{t\in[0,T[}\lmin \bigl(D^2_{xx}
u_t \bigr)> 0 \Bigr)=1.
\end{equation}
Then, the sequence of strategies $\cT_0=\{\cT^n(0):n\geq0\}$ defined
in \eqref{eq1} with $\pa=0$ is admissible and asymptotically optimal,
\[
\lim_{ n\rightarrow+\infty} N_T^{n}\bigl\langle
Z^n\bigr\rangle_T= \biggl( \int_0^T
\tr(X_t) \dt \biggr)^2 \qquad\mbox{a.s.}
\]
\end{theorem}
For the proof, see Section~\ref{subsection:proof:th2:bis}. The extra
assumption \eqref{eq:th3:ass} is satisfied in dimension one for
call/put option in Black--Scholes model only if the hedging time
horizon is strictly smaller than the option maturity. But it is not
satisfied in digital call/put option. This discussion can be extented
to higher multidimensional situations.

%
\begin{remark}
In the one dimensional case, we have
\[
X_t=\frac{1} {\sqrt6}\sigma^2_t\bigl|D_{xx}^2
u_t\bigr|,\qquad \Lambda_t=\frac{1} {
\sqrt6}\bigl|D_{xx}^2
u_t\bigr|,
\]
and the $\mu$-optimal stopping times read
\[
\tau^n_{i}=\inf \biggl\{t\geq\tau^n_{i-1}
\dvtx|S_t - S_{\tau^n_{i-1}}| > \frac
{\en
}{\sqrt{|D_{xx}^2 u_{\tau^n_{i-1}}|/\sqrt6+ \pa\appi(|D_{xx}^2
u_{\tau^n
_{i-1}}|/\sqrt6)}} \biggr\}\wedge T.
\]
For $|D_{xx}^2 u_t|$ bounded from below, we can take $\pa=0$ and the
optimal strategy coincides with that of \cite{F11}, Theorem C.

The
threshold $\pa\neq0$ ensures that the hedging rebalancing occurs often
enough, even if $\Lambda_{t}\neq0$ for some time $t$: this
interpretation is also valid in the multidimensional case.
\end{remark}

\subsection{Proof of Theorem \texorpdfstring{\protect\ref{th1}}{3.1}}
It is split into several steps.

\textit{Step \textup{1:} Quadratic variation decomposition}.
We start from the hedging error~\eqref{eq:zs}. A natural idea consists
in writing a Taylor expansion (regarding the $S$ variable only) and
showing that the residual terms converge to $0$ fast enough as we could expect,
%
\begin{equation}
\label{eq:zn1} Z^n_{s}=\int_0^{s}
\bigl(D^2_{xx} u_{\varphi(t)} \Delta S_t
\bigr)\cdot\dS_t + R^{n}_{s},
\end{equation}
where
%
\begin{equation}
\label{eq:R}R^{n}_{s}:= \int_0^{s}
\bigl(D_x u_t- D_x u_{\varphi(t)}-
D^2_{xx} u_{\varphi(t)} \Delta S_t \bigr)
\cdot\dS_t,\qquad s\leq T.
\end{equation}
Then passing to quadratic variation, we obtain
\[
\bigl\langle Z^n\bigr\rangle_T= \int
_0^T \Delta S_t^*D^2_{xx}
u_{\varphi(t)} \,\mathrm{ d}\langle S \rangle_t D^2_{xx}
u_{\varphi(t)} \Delta S_t +e_{1,T}^n,
\]
where
%
\begin{equation}
\label{eq:e1}e_{1,T}^n:= \bigl\langle R^{n}
\bigr\rangle_T+2 \biggl\langle\int_0^{\cdot}
\bigl(D^2_{xx} u_{\varphi(t)} \Delta S_t
\bigr)\cdot\dS_t,R^{n}_{\cdot} \biggr\rangle_T.
\end{equation}
Now, we wish an expression involving only the Brownian motion for ease
of mathematical analysis: hence we replace $\Delta S_t$ by $\sigma
_{\varphi
(t)} \Delta B_t$ and $\,\mathrm{ d}\langle S\rangle_t$ by $\sigma
_{\varphi
(t)}\sigma
_{\varphi(t)}^*\dt$, leading to
%
\begin{eqnarray}
\label{eq:e2}
\nonumber
\bigl\langle Z^n\bigr\rangle_T&=&
\int_0^T \Delta B_t^* \bigl(
\sigma_{\varphi(t)}^*D^2_{xx} u_{\varphi(t)}
\sigma_{\varphi(t)} \bigr)^2\Delta B_t
\dt+e_{1,T}^n+e_{2,T}^n,
\\
e_{2,T}^n&:=&\int_0^T
\Delta S_t^*D^2_{xx} u_{\varphi(t)} \Delta
\bigl(\sigma_t\sigma_t^* \bigr) D^2_{xx}
u_{\varphi(t)} \Delta S_t\dt
\nonumber
\\[-8pt]
\\[-8pt]
\nonumber
&&{}+
\int_0^T (\Delta
S_t + \sigma_{\varphi(t)} \Delta B_t
)^*\\
&&\hspace*{28pt}{}\times D^2_{xx} u_{\varphi(t)}\sigma_{\varphi(t)}
\sigma_{\varphi
(t)}^*D^2_{xx} u_{\varphi(t)} \biggl(
\int_{\varphi(t)}^t \Delta \sigma_r
\,\mathrm{ d}B_r \biggr)\dt.\nonumber
\end{eqnarray}
As mentioned before, we seek a \emph{smart} representation of the main
term of $\langle Z^n\rangle_T$ in the form $\sum_{\tau_{i-1}^n<
T}(\Delta B_{\tau^n_i}^*X_{\tau^n_{i-1}}
\Delta B_{\tau^n_i} )^2$ plus a stochastic integral, where $X$ is a
measurable adapted $d\times d$-matrix process which has to be defined.
Instead of directly giving the solution, let us discuss a bit on the
expected properties of $X$. Applying It\^o's formula on each interval
$[\tau^n_{i-1},\tau^n_i]$, we obtain
\begin{eqnarray*}
&&\sum_{\tau_{i-1}^n< T} \bigl(\Delta B_{\tau^n_i}^*X_{\tau^n_{i-1}}
\Delta B_{\tau^n_i} \bigr)^2\\
&&\qquad=  \int_0^T
\Delta B_t^* \bigl(2 \tr(X_{\varphi(t)} ) X_{\varphi
(t)} +
\bigl(X_{\varphi(t)}+X_{\varphi(t)}^* \bigr)^2 \bigr)\Delta
B_t \dt
\\
&&\qquad\quad{} +2\int_0^T \Delta B_t^*X_{\varphi(t)}
\Delta B_t\Delta B_t^* \bigl(X_{\varphi(t)}+X_{\varphi(t)}^*
\bigr) \dB_t,
\end{eqnarray*}
with the tentative identification
%
\begin{equation}
\label{eq:tentative}2 \tr(X_{\varphi(t)} ) X_{\varphi(t)} +
\bigl(X_{\varphi(t)}+X_{\varphi(t)}^* \bigr)^2= \bigl(\sigma
_{\varphi
(t)}^*D^2_{xx} u_{\varphi(t)}\sigma
_{\varphi(t)} \bigr)^2.
\end{equation}
Mainly, two reasons prompt us to impose $X_{\varphi(t)}\in\mathcal
{S}_+^d(\R)$.
\begin{itemize}
\item Gathering the previous identities and anticipating a little bit
on the following, the main contribution in ${N^n_T}\langle Z^n\rangle
_T$ is
\[
{N^n_T}\sum_{\tau_{i-1}^n< T} \bigl(
\Delta B_{\tau^n_i}^*X_{\tau
^n_{i-1}} \Delta B_{\tau^n_i}
\bigr)^2\geq \biggl(\sum_{\tau_{i-1}^n< T}\bigl|\Delta
B_{\tau^n_i}^*X_{\tau^n_{i-1}} \Delta B_{\tau^n
_i}\bigr |
\biggr)^2
\]
using the Cauchy--Schwarz inequality. In general the limit of the above
lower bound is not easy to handle because of the absolute values, but
if the matrix $X_{\varphi(t)}$ is nonnegative-definite, we can remove
them and conclude using a convergence result about discrete quadratic
variations (Proposition \ref{prop3}).
\item Once that we have restricted to nonnegative-definite matrices,
let us prove that the solution to \eqref{eq:tentative} (whenever it
exists) is symmetric. If $\tr(X_{\varphi(t)})=0$, then $X_{\varphi
(t)}=0$ (thus symmetric): indeed, $X_{\varphi(t)}+X_{\varphi(t)}^*$ is
symmetric nonnegative-definite and has a null trace, thus it
is the zero-matrix and consequently $X_{\varphi(t)}=-X_{\varphi
(t)}^*=0$ (since both $X_{\varphi(t)}$ and $X_{\varphi(t)}^*$ are
nonnegative-definite). If $\tr(X_{\varphi(t)})>0$, then taking the
transposition of \eqref{eq:tentative} readily gives $X_{\varphi
(t)}=X_{\varphi(t)}^*$.
\end{itemize}
From Lemma \ref{lem5}, there exists exactly one adapted process $X$
with values in $\S_+^d(\R)$, solution of the equation $2 \tr
(X)X+4X^2=(\sigma^* D^2_{xx} u \sigma)^2$.
In addition, this solution is continuous \emphh{a.s.} because
$C:=\sigma^* D^2_{xx}
u \sigma$ is continuous a.s., and the solution $X$ is continuous as a
function of $C$ on $\S^d$. Gathering the previous identities, we have
established a nice decomposition of the quadratic variation of the
hedging error
%
\begin{eqnarray}
\label{eq:qvZT:bis} \bigl\langle Z^n\bigr\rangle_T&= & \sum
_{\tau_{i-1}^n<
T} \bigl(\Delta B_{\tau^n_i}^*X_{\tau^n_{i-1}}
\Delta B_{\tau^n
_i} \bigr)^2 +e_{1,T}^n+e_{2,T}^n+e_{3,T}^n,
\\
e_{3,T}^n&:=& -4\int_0^T
\Delta B_t^*X_{\varphi(t)} \Delta B_t\Delta
B_t^* X_{\varphi(t)} \dB_t. \label{eq:e3}
\end{eqnarray}

\textit{Step \textup{2:} Lower bound for the renormalized quadratic variation}.
The Cauchy--Schwarz inequality yields that ${N^n_T}\sum_{\tau
_{i-1}^n< T}(\Delta
B_{\tau^n
_i}^*X_{\tau^n_{i-1}} \Delta B_{\tau^n_i} )^2$ is bounded from below by
\begin{eqnarray*}
\biggl( \sum_{\tau_{i-1}^n< T}\bigl|\Delta B_{\tau^n_i}^*X_{\tau
^n_{i-1}}
\Delta B_{\tau^n
_i} \bigr| \biggr)^2&=& \biggl(\sum
_{\tau_{i-1}^n< T}\Delta B_{\tau
^n_i}^*X_{\tau^n_{i-1}} \Delta
B_{\tau^n
_i} \biggr)^2\\
&\cvas &\biggl(\int_0^T
\tr(X_t )\dt \biggr)^2,
\end{eqnarray*}
using that $X$ is a nonnegative-definite matrix process and applying
Proposition \ref{prop3}.

\textit{Step \textup{3:} The renormalized errors ${\en^{-2\rhos}}
e_{1,T}^n$, ${\en^{-2\rhos}} e_{2,T}^n$ and ${\en^{-2\rhos}} e_{3,T}^n$
converge to 0} a.s.
Observe that once these convergences are established, in view of \eqref
{eq:qvZT:bis} and \AN\ we easily complete the proof of Theorem \ref
{th1}.

\emph{Proof of ${\en^{-2\rhos}} e_{1,T}^n\cvas0$.}
We first state an intermediate result which is proved in Appendix
(Section~\ref{proof:lem2}).
%
%
\begin{lemma}
\label{lem2} Assume hypotheses \textup{\Abs}, \textup{\Au}, \textup{\AS}, \textup{\AN}\
and \textup{\AY}\ are in
force. Then
${\en^{2-4\rhos}} \langle R^n \rangle_T\cvas0$ where $R^n$ is
defined in \eqref{eq:R}.
\end{lemma}\eject
Then, starting from \eqref{eq:e1}, applying the Cauchy--Schwarz
inequality to the cross-variation and using \Abs--\Au--\AS, we derive
\begin{eqnarray*}
&&{\en^{-2\rhos}}\bigl |e_{1,T}^n\bigr|
\\
&&\qquad\leq{\en^{-2\rhos}} \bigl\langle R^n \bigr
\rangle_T \\
&&\qquad\quad{}+ 2 \biggl({\en^{-2}} \int_0^T
\Delta S^*_t D^2_{xx} u_{\varphi(t)}
\,\mathrm{ d}\langle S \rangle_t D^2_{xx}
u_{\varphi(t)} \Delta S_t \biggr)^{1/2}\ \bigl({
\en^{2-4\rhos}} \bigl\langle R^n \bigr\rangle_T
\bigr)^{1/2}
\\
&&\qquad \leq{\en^{2(\rhos-1)}\en^{2-4\rhos}} \bigl\langle R^n \bigr
\rangle_T + 2 C_0 \bigl({\en^{2-4\rhos}} \bigl
\langle R^n \bigr\rangle_T \bigr)^{1/2} \cvas0.
\end{eqnarray*}

\emph{Proof of ${\en^{-2\rhos}} e_{2,T}^n\cvas0$.}
We analyze separately the two contributions in~\eqref{eq:e2}.
\begin{longlist}[(1)]
\item[(1)] First, simple computations using \Abs--\Au--\AS\ {and
Corollary \ref{cor2} directly give (for any given $\rho>0$)}
\begin{eqnarray*}
&&{\en^{-2\rhos}} \biggl|\int_0^T \Delta
S_t^*D^2_{xx} u_{\varphi(t)} \Delta \bigl(
\sigma_t\sigma_t^* \bigr) D^2_{xx}
u_{\varphi(t)} \Delta S_t\dt\biggr|\\
&&\qquad \leq C_0 {
\en^{-2\rhos+2} \bigl(\en^{\ths}+\en^{({\ths}/2)(2-\rho
)} \bigr)}.
\end{eqnarray*}
{Since $\rhos< 1+\ths/2$ and $\rho$ can be taken arbitrarily small, we
obtain that the above upper bound converges \emphh{a.s.} to 0.}
\item[(2)] Second, we apply twice Corollary \ref{cor3}{(ii)}, first taking $\theta=0$ and second taking $\theta=\ths$, so that we
obtain, {for any given $\rho>0$,} \emphh{a.s.} for any $n\geq0$,
%
\begin{eqnarray}
&&\supi\supit\bigl|\Delta S_t+\sigma_{\varphi(t)}\Delta B_t\bigr|
\leq C_0\en^{1-\rho
},\label{eq:accrois:proof:th1:1}
\\
&&\supi\supit\biggl|\int_{\varphi(t)}^t \Delta\sigma
_r \,\mathrm{ d}B_r\biggr| \leq C_0 \en
^{1+\ths-\rho},\label{eq:accrois:proof:th1:2}
\\
&&{\en^{-2\rhos}} \biggl|\int_0^T (\Delta
S_t+\sigma_{\varphi
(t)}\Delta B_t)^*\nonumber\\
&&\hspace*{46pt}{}\times D^2_{xx}
u_{\varphi(t)}\sigma_{\varphi(t)}\sigma_{\varphi
(t)}^*D^2_{xx}
u_{\varphi(t)} \biggl(\int_{\varphi(t)}^t \Delta
\sigma _r \,\mathrm{ d}B_r \biggr)\dt\biggr|
\nonumber
\\
&&\qquad\leq{C_0 \en^{2+\ths-2\rhos-2\rho}}.
\nonumber
\end{eqnarray}
\end{longlist}
{Owing to $\rhos< 1+\ths/2$, taking $\rho$ small enough implies the
\emphh{a.s.}
convergence of the latter upper bound to 0. As a result, $\en^{-2\rhos}
e_{2,T}^n\cvas0$}. 

\emph{Proof of ${\en^{-2\rhos}} e_{3,T}^n\cvas0$.} It is
a direct consequence of the following lemma.  
%
%
\begin{lemma}
Assume \textup{\Abs}. Let $\cT=(\cT^n)_{n\geq0}$ be an admissible sequence of
strategies, and let $(H_t)_{0\leq t < T}$ be a continuous adapted
$d\times d$-matrix process such that $\sup_{t\in[0,T)}|H_t|<+\infty$
\emphh{a.s.}
Then for any {$p>\frac{2}{3-2\rhos}$}, the series $\sum_{n\geq0}
|{\en^{-2\rhos}} \int_{0}^{T} \Delta B_t^* H_{\varphi(t)} \Delta
B_t\Delta B_t^*H_{\varphi(t)} \dB_t|^{p}$ converges almost surely.
\end{lemma}
\begin{pf}
Set $\alpha^n_t:=\Delta B_t^* H_{\varphi(t)} \Delta B_t\Delta
B_t^*H_{\varphi(t)}$ and define the scalar continuous local martingale
$M^n_t:={\en^{-2\rhos}}\int_{0}^t \alpha^n_s \dB_s$. In view of
Corollary \ref{cor1}, it is enough to check that $(\langle M^n \rangle
^{p/2}_T)_{n\geq0}$ defines the terms of an \emphh{a.s.} convergent
series. An
application of Corollary \ref{cor3}(ii) with {$\rho=\frac
{(3-2\rhos)p-2}{3p}>0$} and $\theta=0$ gives $\supi\supit
|\Delta B_t|<C_0 \en^{1-\rho}$ and therefore
\begin{eqnarray*}
\bigl\langle M^n \bigr\rangle_T^{p/2}&=&
\en^{-2p\rhos} \biggl(\int_0^T\bigl |
\alpha^n_t\bigr|^2 \dt \biggr)^{p/2}\\
&\leq&
C_0 \en^{-2p\rhos}\supi\supit|\Delta B_t|^{3p}
\leq C_0\en^{2} \qquad\mbox{a.s.}
\end{eqnarray*}
We are finished.
\end{pf}

\subsection{Proof of Theorem \texorpdfstring{\protect\ref{th2}}{3.2}}
\label{subsection:proof:th2}
We first check the admissibility of $\cT_\pa$, by applying Proposition
\ref{prop4}. Indeed, owing to \Au\ and \Abse, $(\Lambda_t)_{0\leq t
<T}$ is a continuous adapted
nonnegative-definite $d\times d$-matrix process with\break  $\sup_{0\leq
t<T}|\Lambda_t|<+\infty$ \emphh{a.s.} The same properties clearly
hold for
$(\Lambda^\pa_t)_{0\leq t <T}$. In~addition, $\lmin(\Lambda^\pa
_t)\geq
\pa/2>0$ and $\sup_{0\leq t <T}\lmax(\Lambda^\pa_t)\leq\pa+\break  \sup_{0\leq t <T}\lmax(\Lambda_t)<+\infty$ \emphh{a.s.} Therefore, $\cT
_\pa$ is
admissible and in addition $\supn\en^2 {N^n_T}<+\infty$ \emphh{a.s.}
Hence it
allows us to re-use the computations of the proof of Theorem \ref{th1}
in the case $\rhos=1$.

Now let us show the $\mu$-optimality.
Writing ${N^n_T}=1+\sum_{1\leq i \leq{N^n_T}-1} 1$, we point out
%
\begin{eqnarray}\label{eq:opt:1}
\en^{2} {N^n_T}&=& \en^{2} + \sum
_{1\leq i \leq{N^n_T}-1} \Delta S_{\tau^n_{i}}^*
\Lambda^\pa_{\tau^n_{i-1}} \Delta S_{\tau^n_{i}}
\nonumber
\\
&=&\en^2 - \Delta S_T^* \Lambda^\pa_{\tau^n_{{N^n_T}-1}}
\Delta S_T + \sum_{\tau_{i-1}^n< T}\Delta
S_{\tau^n_{i}}^* \Lambda ^\pa_{\tau^n_{i-1}} \Delta
S_{\tau^n_{i}} \\
&\cvas&\int_0^T \tr \bigl(
\Lambda^\pa_t \sigma_t \sigma_t^*
\bigr)\dt\nonumber
\end{eqnarray}
using the convergence of Proposition \ref{prop3}.
On the other hand, starting from the decomposition \eqref{eq:qvZT:bis}
of the hedging error quadratic variation, we write
%
\begin{eqnarray}\label{eq:qvZT:ter:2}
\bigl\langle Z^n\bigr\rangle_T&=& \sum
_{1\leq i \leq{N^n_T}-1} \bigl(\Delta S_{\tau^n_i}^* \Lambda^\pa_{\tau^n_{i-1}}
\Delta S_{\tau^n_i} \bigr)^2 +e_{1,T}^n+e_{2,T}^n+e_{3,T}^n\nonumber\\
&&\hspace*{40pt}{}+e_{4,T}^n+e_{5,T}^n
+e_{6,T}^n,
\nonumber\\
e_{4,T}^n&:=&\sum_{\tau_{i-1}^n< T} \bigl(
\Delta B_{\tau
^n_i}^*X_{\tau^n_{i-1}} \Delta B_{\tau^n_i}
\bigr)^2- \bigl(\Delta S_{\tau^n_i}^* \Lambda_{\tau^n_{i-1}}
\Delta S_{\tau^n
_i} \bigr)^2,
\\
e_{5,T}^n&:=&\sum_{\tau_{i-1}^n< T} \bigl(
\Delta S_{\tau^n_i}^* \Lambda_{\tau^n_{i-1}} \Delta S_{\tau^n_i}
\bigr)^2- \bigl( \Delta S_{\tau^n_i}^*\Lambda^\pa_{\tau^n_{i-1}}
\Delta S_{\tau^n_i} \bigr)^2,
\nonumber\\
e_{6,T}^n&:=& \bigl(\Delta S_{T}^*
\Lambda^\mu_{\tau^n_{{N^n_T}-1}} \Delta S_{T}
\bigr)^2.
\nonumber
\end{eqnarray}
In view of the definition of the strategy $\cT^n_\pa$, \eqref
{eq:qvZT:ter:2} becomes
%
\begin{equation}
\en^{-2}\bigl\langle Z^n\bigr\rangle_T=
\sum_{1\leq i \leq{N^n_T}-1} \Delta S_{\tau^n
_i}^*
\Lambda^\pa_{\tau^n_{i-1}} \Delta S_{\tau^n_i} +
\en^{-2}\sum_{j=1}^6
e_{j,T}^n\label{eq:qvZT:ter:3}.
\end{equation}
Similarly to \eqref{eq:opt:1}, we show that $\sum_{1\leq i \leq{N^n_T}-1}
\Delta S_{\tau^n_i}^*\Lambda^\pa_{\tau^n_{i-1}} \Delta S_{\tau
^n_i} \cvas\break
\int_0^T \tr(\Lambda^\pa_t \sigma_t\sigma_t^* )\dt$.
Furthermore we have already established (see step 3 of proof of Theorem
\ref{th1}) that $\en^{-2}e_{j,T}^n\cvas0$ for $j=1,2,3$ {(remind that
we can take $\rhos=1$); }
the case $j=6$ is also fulfilled because $0\leq e_{6,T}^n\leq\en^4$.

To analyze $e_{4,T}^n$, set $D_{B,i}:=\si\Delta B_{\tau^n_i}$ and
$D_{S,i}:=\Delta S_{\tau^n_i}$, write $X_{\tau^n_{i-1}}=\si^*\li\si
$ and
\begin{eqnarray*}
&&\bigl(\Delta B_{\tau^n_i}^*X_{\tau^n_{i-1}} \Delta B_{\tau^n_i}
\bigr)^2- \bigl(\Delta S_{\tau^n_i}^*\Lambda_{\tau^n_{i-1}}
\Delta S_{\tau^n_i} \bigr)^2\\
&&\qquad= \bigl(\,\dbi^*\li\dbi
\bigr)^2- \bigl(\dsi^*\li\dsi \bigr)^2
\\
&&\qquad= \bigl(\dbi^*\li\dbi-\dsi^*\li\dsi \bigr) \bigl(\dbi^*\li \dbi+\dsi^*\li\dsi
\bigr)
\\
&&\qquad= (\dbi+\dsi)^*\li(\dbi-\dsi) \bigl(\dbi^*\li\dbi+\dsi^*\li \dsi \bigr).
\end{eqnarray*}
Then we deduce that $\en^{-2}|e_{4,T}^n|$ is bounded by
\begin{eqnarray*}
&&\en^{-2}{N^n_T}\supi\sup
_{\tau^n_{i-1}\leq t \leq\tau
^n_i]}|\Lambda _{\varphi
(t)}|^2|\Delta
S_{t}+\sigma_{\varphi(t)} \Delta B_t| \biggl|\int
_{\varphi
(t)}^t \Delta\sigma _s
\dB_s \biggr| \\
&&\hspace*{109pt}{}\times\bigl(|\Delta S_{t}|^2+|
\sigma_{\varphi
(t)}\Delta B_t|^2 \bigr)
\\
&&\qquad\leq C_0\en^{-2}\en^{-2}\en^{1-\rho}
\en^{(1+\ths-\rho)}\en^{2(1-\rho
)}=C_0\en^{\ths/5}\cvas0,
\end{eqnarray*}
where we have used \AN\ {(with $\rhos=1$)}, and estimates (\ref
{eq:accrois:proof:th1:1})--(\ref{eq:accrois:proof:th1:2}) with $\rho
=\ths
/5$ (which are available for any sequence of admissible strategies).
This proves $\en^{-2}e_{4,T}^n\cvas0$.

Finally regarding $e_{5,T}^n$, recalling that the matrix $\li$ is
nonnegative-definite, we obtain that $|\en^{-2}e_{5,T}^n|$ is bounded by
\begin{eqnarray*}
&& \en^{-2}\sum_{\tau_{i-1}^n< T}\bigl|\Delta
S_{\tau^n_i}^* \Lambda _{\tau^n_{i-1}} \Delta S_{\tau^n
_i} -\Delta
S_{\tau^n_i}^*\Lambda^\pa_{\tau^n_{i-1}} \Delta
S_{\tau^n_i}\bigr |\\
&&\hspace*{44pt}{}\times \bigl(\Delta S_{\tau^n_i}^* \Lambda_{\tau^n_{i-1}}
\Delta S_{\tau^n_i} +\Delta S_{\tau^n_i}^*\Lambda^\pa_{\tau^n_{i-1}}
\Delta S_{\tau^n_i} \bigr)
\\
&&\qquad\leq\sum_{\tau_{i-1}^n< T}\mu\appi \bigl(\lmin(
\Lambda_{\tau
^n_{i-1}}) \bigr) |\Delta S_{\tau^n_i}|^2 \bigl[2
\en^{-2} \Delta S_{\tau^n_i}^*\Lambda^\pa_{\tau^n_{i-1}}
\Delta S_{\tau^n
_i} \bigr]
\\
&&\qquad\leq2 \mu\sum_{\tau_{i-1}^n< T}\appi \bigl(\lmin(
\Lambda_{\tau
^n_{i-1}}) \bigr) |\Delta S_{\tau^n_i}|^2,
\end{eqnarray*}
where we have used the definition of $\cT_\pa$ at the last inequality.
Thus Proposition~\ref{prop3} yields
\[
\limsup_{n\rightarrow+\infty}\bigl|\en^{-2}e_{5,T}^n\bigr|
\leq2\mu\int_0^T \appi \bigl(\lmin(
\Lambda_{t}) \bigr)\tr \bigl(\sigma_t\sigma_t^*
\bigr)\dt\qquad \mbox{a.s.}
\]
Let us summarize: setting $L_T:=\int_0^T \tr(\Lambda_t \sigma
_t\sigma_t^* )\dt=\int_0^T \tr(X_t )\dt$ and $L^\pa
_T:=\int_0^T \appi(\lmin(\Lambda_{t})) \tr(\sigma
_t\sigma
_t^*)\dt$ so that
$\int_0^T \tr(\Lambda^\pa_t \sigma_t\sigma_t^* )\dt
=L_T+\mu
L^\pa_T$,
we have shown
\begin{eqnarray*}
&&\en^{2} {N^n_T}\cvas L_T+\mu
L^\pa_T,\qquad \limsup_{n\rightarrow+\infty} \bigl|
\en^{-2}\bigl\langle Z^n\bigr\rangle_T-
\bigl(L_T+ \mu L^\pa_T \bigr) \bigr|\leq2\pa
L^\pa_T \qquad\mbox{a.s.},
\\
&&\limsup_{n\rightarrow+\infty} \bigl|{N^n_T}\bigl\langle
Z^n\bigr\rangle_T-(L_T)^2\bigr|
\\
&&\qquad\leq\limsup_{n\rightarrow+\infty} \bigr|\en^{-2} \bigl\langle
Z^n\bigr\rangle_T-L_T\bigr| \limsup
_{n\rightarrow+\infty} \en^{2}{N^n_T}+L_T
\limsup_{n\rightarrow+\infty} \bigl|\en^2{N^n_T}-L_T\bigr|
\\
&&\qquad\leq3 \mu L^\pa_T \bigl(L_T+\mu
L^\pa_T \bigr)+L_T \mu L^\pa_T=
\mu L^\pa_T \bigl(4L_T+3\mu
L^\pa_T \bigr)\qquad\mbox{a.s.}
\end{eqnarray*}
Theorem \ref{th2} is proved.

\subsection{Proof of Theorem \texorpdfstring{\protect\ref{th2:bis}}{3.3}}
\label{subsection:proof:th2:bis}
Here, arguments are simpler in all steps of the proof of Section~\ref{subsection:proof:th2}, so we shall skip details; the admissibility of
the strategy comes readily from the ad hoc assumption \eqref
{eq:th3:ass} and Proposition \ref{prop4}; the optimality follows as
before from
\[
\en^{2} {N^n_T}= \en^{2} + \sum
_{1\leq i \leq{N^n_T}-1} \Delta S_{\tau^n_{i}}^*
\Lambda_{\tau^n_{i-1}} \Delta S_{\tau^n_{i}}
\cvas\int
_0^T \tr(X_t ) \dt,
\]
and from [setting $\bar e_{6,T}^n:=(\Delta S_{T}^* \Lambda_{\tau^n_{{N^n_T}
-1}} \Delta S_{T})^2$]
\begin{eqnarray*}
\en^{-2}\bigl\langle Z^n\bigr\rangle_T&=&
\en^{-2}\sum_{1\leq i \leq{N^n_T}-1} \bigl(\Delta
S_{\tau^n
_i}^*\Lambda_{\tau^n_{i-1}} \Delta S_{\tau^n_i}
\bigr)^2 +\en^{-2}\sum_{j=1}^4
e_{j,T}^n+\en^{-2}\bar e_{6,T}^n
\\
&\cvas&\int_0^T \tr(X_t )\dt
\end{eqnarray*}
with the help of the convergence results already obtained. Theorem \ref
{th2:bis} is proved. 

\section{Numerical experiments}
\label{sec4}
\subsection{Algorithm for the optimal stopping times}
From the previous section (Theorem \ref{th2}), the $\mu$-optimal
stopping times $(\mu>0)$ are iteratively given by $\tau^n_{0}:=0$ and
\[
\tau^n_{i}:=\inf \bigl\{t\geq\tau^n_{i-1}
\dvtx(S_t - S_{\tau^n_{i-1}})^* \Lambda^\pa_{\tau^n_{i-1}}
(S_t - S_{\tau^n_{i-1}}) \geq\en^2 \bigr\}\land T,
\]
where for any $t$, $\Lambda^\pa_t:=\Lambda_t + \pa\appi(\lmin
(\Lambda
_t))I_d$, $\Lambda_t:=(\sigma^{-1}_t)^*X_t \sigma^{-1}_t$ and
$X_t$ solves
\eqref{eq5} with $c_t=\sigma_t^* D^2_{xx} u_t \sigma_t$. Thus,
$\tau^n_i$ is the
first hitting time of an ellipsoid centered at $S_{\tau^n_{i-1}}$ with
principal axes equal to the orthogonal eigenvectors of the symmetric
positive-definite matrix $\Lambda^\pa_{\tau^n_{i-1}}$ (or equivalently
those of $\Lambda_{\tau^n_{i-1}}$).
We briefly recall (see Section~\ref{subsection:proof:lem5}) the main
steps to compute the matrix $X_{\tau^n_{i-1}}$ ($i\geq1$) from which we
derive $\Lambda_{\tau^n_{i-1}}$ and $\Lambda^\mu_{\tau^n_{i-1}}$:
\begin{longlist}[(1)]
\item[(1)] Diagonalize the symmetric matrix $c_{\tau^n_{i-1}}=\sigma
_{\tau^n
_{i-1}}^* D^2_{xx} u_{\tau^n_{i-1}} \sigma_{\tau^n_{i-1}}:=P_{\tau^n_{i-1}} \times \operatorname{Diag} (\lambda_j(c_{\tau^n
_{i-1}})\dvtx1\leq j \leq d ) P^*_{\tau^n_{i-1}}$, where $P_{\tau
^n_{i-1}}$ is
an orthogonal matrix.
\item[(2)] Find the zero $y_{\tau^n_{i-1}}\in\R^+$ of the increasing
function $ y\mapsto(4+d)y-\sum_{j=1}^d \sqrt{y^2+4\lambda
_j^2(c_{\tau^n_{i-1}})}$. This root lies in the interval
$
[0,d|\lambda(c_{\tau^n_{i-1}})|/\sqrt{4+2d} ]$; see the proof of Lemma
\ref{lem5}.
\item[(3)] From \eqref{eq:sol:xc}, we obtain
\[
X_{\tau^n_{i-1}}=P_{\tau^n_{i-1}} \operatorname{Diag} \biggl( \frac
{-y_{\tau^n_{i-1}}+\sqrt{y_{\tau^n_{i-1}}^2+4\lambda^2_j(c_{\tau^n
_{i-1}})}}{4}
\dvtx1\leq j \leq d \biggr) P^*_{\tau^n_{i-1}}.
\]
\end{longlist}
Last, we mention that even if $\Lambda^\pa_{\tau^n_{i-1}}$ is tractable,
the exact simulation of $\tau^n_{i}$ is in generally impossible, and
approximations are required; see \cite{gobe:meno:10} and references therein.
\subsection{Numerical tests}
This section is dedicated to an application of Theorem~\ref{th2} to the
case of an exchange binary option $g(S_T)=\1_{S^1_T\geq S^2_T}$. This
example is relevant in our study (and improves the setting of \cite
{F11}) because this is a simple \emph{bi-dimensional nonconvex}
function, for which the value function $u$ and its sensitivities are
available in the Black--Scholes model
\[
\mathrm{ d} %
\pmatrix{S^1_t\vspace*{2pt}
\cr
S^2_t } %
= %
\pmatrix{\sigma
_1S^1_t&0\vspace*{2pt}
\cr
\rho\sigma
_2S^2_t&\sqrt{1-\rho^2}\sigma
_2S^2_t } %
\mathrm{ d} %
\pmatrix{B^1_t\vspace*{2pt}
\cr
B^2_t
},
\]
where $(B^1,B^2)$ are two independent Brownian motions. The model
parameters are set to $S^1_0=100$, $S^2_0=100$, $\sigma_1=0.3$,
$\sigma_2=0.4$,
$\rho=0.5$ and $T=1$.

We take $\en=0.05$. In our different tests, we have not observed a
significant difference by taking $\mu=0$ or $\mu$ small; hence, we only
report the values for $\mu=0$. We generate 1000 experiments $\omega$,
independently. To compute the hitting times for each $\omega$, we use a
thin uniform time
mesh $\pi_{\bar n}=(i T/\bar n)_{0\leq i\leq\bar n}$ ($\bar n=50\mbox{,}000$
in our tests): we draw $S^1(\omega)$ and $S^2(\omega)$ along $\pi
_{\bar n}$ and compute (with the help of the previous algorithm) the
hitting times $\tau^n_i(\omega)=\inf\{t\in\pi_{\bar n}\,\cap\,]\tau^n
_{i-1}(\omega),T]\dvtx[(S_t - S_{\tau^n_{i-1}})^* \Lambda^\pa
_{\tau^n_{i-1}} (S_t
- S_{\tau^n_{i-1}})](\omega) \geq\en^2 \}\land T$; at the end of the
process, we get the number ${N^n_T}(\omega)$ of discrete times. The mesh
$\pi_{\bar n}$ is also used to compute subsequent quadratic variations
and time integrals.

We compare $\omega$ by $\omega$ the above strategy with that based on
the \emph{uniform} mesh $\pi_{{N^n_T}(\omega)}$ and with that based
on the
so-called \emph{fractional} mesh\footnote{According to \cite{GG11}, the
fractional smoothness of $g(S_T)$ is $\frac{ 1 }{2 }$; thus, when ${N^n_T}
(\omega)$ is deterministic, this choice of fractional mesh yields that
$\E(\langle Z^n\rangle_T)$ is of order 1 w.r.t. the inverse of the
number of times, instead of order $\frac{ 1 }{2 }$ with the uniform mesh.}
$ (T [1-(1-i/{N^n_T}(\omega))^2 ] )_{1\leq
i\leq{N^n_T}
(\omega)}$:
this comparison looks quite fair from a practitioner point of view
since he is allowed to rebalance the hedging portfolio the same number
of times. The use of the optimal stochastic grid is slightly more
demanding since it requires the computations of more Greeks than only
the Delta (because of the matrix $\Lambda^\mu$); however, these
sensitivities are widely available in any trading system, which makes
this higher complexity likely negligible in view of the benefit of
optimal times.

%
\begin{figure}

\includegraphics{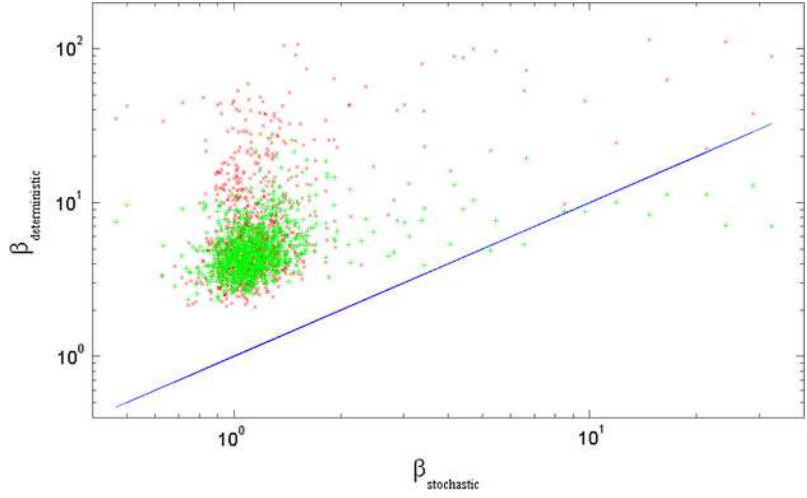}

\caption{``{$\times$},'' ``{$+$}'' and the
blue line correspond respectively to ``{$(\beta
_{\mathrm
{stochastic}},\beta_{\mathrm{uniform}})$},'' ``{$(\beta
_{\mathrm{stochastic}},\beta_{\mathrm{fractional}})$}'' and the identity
function.}
\label{fig2}
\end{figure}

We define $\beta_{\mathrm{stochastic}}(\omega)$, $\beta_{\mathrm
{uniform}}(\omega)$, $\beta_{\mathrm{fractional}}(\omega)$ where we
compute\break
$\beta_{\cdot}(\omega):=\frac{{N^n_T}\langle Z^n\rangle_T}{(\int_0^T \tr(X_t
)\dt)^2}(\omega)$ according to each of these three strategies: in view
of Theorem~\ref{th2}, this ratio is asymptotically greater than 1 and
adimensional; moreover, the closer to 1 the ratio, the better the
strategy.

\emph{Results.}
Figure~\ref{fig2} displays, for each $\omega$, the couples
\[
\bigl(\beta_{\mathrm{stochastic}}(\omega),\beta_{\mathrm
{uniform}}(\omega) \bigr) \quad\mbox{and}\quad \bigl(\beta_{\mathrm{stochastic}}(\omega),\beta_{\mathrm
{fractional}}(\omega) \bigr).
\]
Most of the times, the points are above the diagonal, showing that the
$\mu$-optimal strategy lessens the quadratic variation $\omega$-wise
(remind that the strategies have got the same number of discrete times
${N^n_T}$), compared to the quadratic variation worked out over the
deterministic time mesh.\ In addition, $\beta_{\mathrm{stochastic}}$ is
concentrated around $1$, which means a convergence of ${N^n_T}\langle
Z^n\rangle_T$ toward
the lower bound $(\int_0^T \tr(X_t)\dt)^2$.

\begin{figure}

\includegraphics{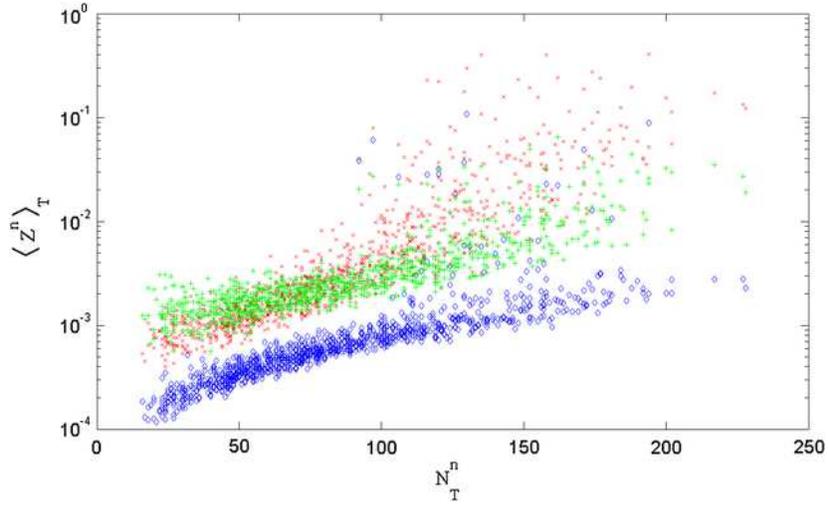}

\caption{``{$\times$},'' ``{$+$}'' and
correspond respectively to
``{$\langle Z^n\rangle_{T,\mathrm{uniform}}$},''
``{$\langle Z^n\rangle_{T,\mathrm{fractional}}$}'' and
``{$\langle Z^n\rangle_{T,\mathrm{stochastic}}$}.''}
\label{fig1}
\end{figure}

Figure~\ref{fig1} displays $\langle Z^n\rangle_T$ as a function of
${N^n_T}$ for
the three strategies and for different $\omega$: here again, we observe
that the $\mu$-optimal strategy outperforms deterministic strategies.\vadjust{\goodbreak}
%
%

\begin{appendix}
\section*{Appendix}
\subsection{Proof of Proposition \texorpdfstring{\protect\ref{prop1}}{2.1}}
\label{appendix:prop1}

Let us prove {(i)}, assuming only \AS. For \mbox{$p=0$},
this is trivial.

Now consider the case $p>0$. Since $\sigma_t$ is
nonzero for any $t$ and continuous, $C_E:=\inf_{t\in[0,T]} (\sum_{j=1}^d e_j\cdot\sigma_t\sigma_t^* e_j)>0$ a.s., where $e_j$ is the $j$th
element of the canonical basis in $\R^d$. Therefore, \emphh{a.s.} for
any $0\leq
s\leq t \leq T$ we have
%
\begin{eqnarray}
\label{eq:CE:1} 0&\leq& t-s \leq C_E^{-1} \int
_{s}^{t} \Biggl(\sum_{j=1}^d
e_j\cdot\sigma_r\sigma_r^* e_j
\Biggr) \,\mathrm{d}r=C_E^{-1 } \sum_{j=1}^d
\bigl[ \bigl\langle S^j \bigr\rangle_{t}- \bigl\langle
S^j \bigr\rangle_{s} \bigr]
\nonumber
\\[-8pt]
\\[-8pt]
\nonumber
&=&C_E^{-1 }\sum_{j=1}^d
\biggl[ \bigl(S^j_t-S^j_s
\bigr)^2 - 2\int_s^t
\bigl(S^j_r-S^j_s \bigr)
\dS^j_r \biggr],
\end{eqnarray}
applying the It\^o formula at the last equality.
Take $s=\tau^n_{i-1}$, $t=\tau^n_{i}$ and use \AS
%
\begin{eqnarray}
\label{eq:CE:2} \Delta\tau^n_i&\leq&
C_E^{-1 } \Biggl(C_0\en^2 +2\sum
_{j=1}^d \biggl|\int_{\tau^n
_{i-1}}^{\tau^n_i}
\Delta S^j_r \dS^j_r \biggr|
\Biggr)
\nonumber
\\[-8pt]
\\[-8pt]
\nonumber
&\leq& C_E^{-1 } \Biggl(C_0\en^2
+4\sum_{j=1}^d \supt\biggl|\int
_{0}^{t} \Delta S^j_r
\dS^j_r \biggr| \Biggr).
\end{eqnarray}
Now for $j=1,\ldots, d$, set $M^{j,n}_{t}:=\en^{2/p-1} \int_{0}^{t}
\Delta S^j_r \dS^j_r$ (recalling that $p>0$). Then
\[
\sum_{n\geq0} \bigl\langle M^{j,n} \bigr
\rangle_T^{p/2}= \sum_{n\geq0}
\en^{2-p} \biggl(\int_0^T \bigl| \Delta
S^j_t\bigr|^2\,\mathrm{ d} \bigl\langle
S^j \bigr\rangle_t \biggr)^{p/2}\leq
C_0 \sum_{n\geq0}\en^{2}<+
\infty\qquad\mbox{a.s.}
\]
Thus owing to Corollary \ref{cor1} the terms $ (\sup_{0\leq t\leq
T} |M^{j,n}_{t} |^p )_{n\geq0}$ define an \emphh{a.s.}
convergent series.
Combining this with \eqref{eq:CE:2}, we finally derive
\begin{eqnarray*}
&&\sum_{n\geq0} \Bigl[\en^{2/p-1}\sup
_{1\leq i \leq{N^n_T}}\bigl|\Delta \tau^n_i\bigr|
\Bigr]^p
\\
&&\qquad\leq\co \Biggl(\sum_{n\geq0} \bigl[\en^{2/p-1}
\en^2 \bigr]^{p}+\sum_{j=1}^d
\sum_{n\geq0} \sup_{0\leq t\leq T}
\bigl|M^{j,n}_{t} \bigr|^p \Biggr)<+\infty \qquad\mbox{a.s.}
\end{eqnarray*}

It remains to justify {(ii)}.
For $p=0$, the result directly follows from \AN\ and the inequality
\eqref{eq:sum:en}. Now
take $p>0$, and set
\begin{eqnarray*}
U_t^n&:=&{\en^{-2(p-1)+2\rhos}}\somt\Biggl|\sum
_{j=1}^d\Delta \bigl\langle S^j \bigr
\rangle_{\tau^n_{i}\wedge t} \Biggr|^p,
\\
V_t^n&:=&{\en^{-2(p-1)+2\rhos}}\somt\ \sup
_{s\in(\tau^n_{i-1},\tau^n_i \wedge t]}\llvert\Delta S_{s}\rrvert^{2p}.
\end{eqnarray*}
If $\sum_{n\geq0}U_T^n \cvas$, \eqref{eq:CE:1} immediately yields that
$\sum_{n\geq0}{\en^{-2(p-1)+2\rhos}}\times\break  \sum_{\tau_{i-1}^n< T}(\Delta
\tau^n_{i})^p \cvas$.
Thus, it is sufficient to show $\sum_{n\geq0}U_t^n \cvas$, for any
$t\in
[0,T]$, and this is achieved by an application of Lemma \ref{lem1:bis}.
The sequences of processes $(U^n)_{n\geq0}$ and $(V^n)_{n\geq0}$ are
in $\C_0^+$.
Then $V^n$ is nondecreasing, and using \AS--\AN
\[
\sum_{n\geq0}V_T^n \leq
C_0\sum_{n\geq0}{ \en^{-2(p-1)+2\rhos}}
\en^{2p}{N^n_T}\leq\co\sum
_{n\geq0}\en ^2<+\infty \qquad\mbox{a.s.}
\]
Then we deduce
that items {(i$'$)} and {(ii$'$)} of Lemma \ref{lem1:bis} are
fulfilled. It remains to check the relation of domination [item {(iii$'$)}].
Let $k\in\bN$. On the set $\{\tau^n_{i-1}<t \wedge\st\}
$, from
the multidimensional BDG inequality in a conditional version, we have
%
\begin{equation}
\label{eq:prop1:bdg}  \E \Biggl( \Biggl|\sum_{j=1}^d
\Delta \bigl\langle S^j \bigr\rangle_{\tau^n
_{i}\wedge
t\wedge\st}\Biggr |^p\Big|
\F_{\tau^n_{i-1}} \Biggr)\leq c_{p} \E \Bigl(\sup
_{\tau^n_{i-1}<
s\leq\tau^n_{i}\wedge t\wedge\st} \llvert\Delta S_{s}\rrvert^{2p}\big|
\F_{\tau^n_{i-1}} \Bigr).\hspace*{-35pt}
\end{equation}
Then, it follows
\begin{eqnarray*}
\E \bigl[ U_{t \wedge\st}^n \bigr] &=&{\en^{-2(p-1)+2\rhos}}\sum
_{i=1}^{+\infty} \E \Biggl(1_{\tau^n_{i-1}<
t\wedge\st}\E
\Biggl[ \Biggl|\sum_{j=1}^d\Delta\langle S
\rangle_{\tau^n
_{i}\wedge t\wedge\st} \Biggr|^p\Big|\F_{\tau^n_{i-1}} \Biggr] \Biggr)
\\
&\leq& c_p \E \bigl[ V_{t \wedge\st}^n \bigr].
\end{eqnarray*}
The proof is complete.

\subsection{Proof of Proposition \texorpdfstring{\protect\ref{prop2}}{2.2}}
\label{appendix:prop2}
Let $p>0$.
Let $\delta$ be the parameter standing for $\frac{1}2$ under \AS\ and
$1$ under \AS--\AN. Set
\begin{eqnarray*}
U_t^n&:=&{\en^{-2\delta(({p(\theta+1)}/2)-2(1-\delta)
)+2+2\rhos
(2\delta-1)}}\somt\sup
_{\tau^n_{i-1} \leq s\leq\tau_i^n\wedge t} \bigl|\Delta M^n_t\bigr|^p,
\\
V_t^n&:=&{\en^{-2\delta((p(\theta+1)/2)-2(1-\delta)
)+2+2\rhos(2\delta-1)}} \somt\biggl|\int
_{\tau^n_{i-1}}^{\tau^n_{i}\wedge t} \alpha^n_r
\,\mathrm{ d}r \biggr|^{p/2}.
\end{eqnarray*}
Observe that the announced result reads as $\sum_{n\geq0}U^n_T\cvas
$. To
prove this convergence, it is enough to establish that $\sum_{n\geq0}
V^n_T\cvas$. Indeed, following the arguments of the proof of
Proposition \ref{prop1}{(ii)}, we can apply Lemma \ref{lem1:bis}
since $(U^n)_{n\geq0}$ and $(V^n)_{n\geq0}$ are two sequences of
continuous adapted processes and:
\begin{longlist}[(ii$'$)]
\item[(i$'$)] $V^n$ is nondecreasing on $[0,T]$ \emphh{a.s.};
\item[(iii$'$)] the domination is satisfied thanks to the BDG
inequalities, similarly to~\eqref{eq:prop1:bdg}.
\end{longlist}
Now to prove (ii$'$), that is, $\sum_{n\geq0}V^n_T\cvas$, write
\begin{eqnarray*}
\sum_{n\geq0}V_T^n &\leq&\sum
_{n\geq0}{\en^{-2\delta(
({p(\theta
+1)}/2)-2(1-\delta) )+2+2\rhos(2\delta-1)}}\\
&& \hspace*{15pt}{}\times\sum
_{\tau_{i-1}^n<
T}\bigl|\co \bigl(\en^{2\theta}+\bigl(\Delta
\tau^n_i\bigr)^\theta \bigr)\Delta
\tau^n_i \bigr|^{p/2} \qquad\mbox{a.s.}
\end{eqnarray*}
First, consider the case \AS\ and set $D_n^{(q)}:=\supi(\Delta\tau^n
_{i})^q$ for $q\geq0$: Proposition \ref{prop1}{(i)} yields
${\cal
D}^{(q)}:=\sum_{n\geq0}\en^{-(q-2)} D_n^{(q)} <+\infty$ \emphh
{a.s.} Using
$p\geq2$,
it readily follows that
\begin{eqnarray*}
\sum_{n\geq0}V_T^n &\leq&\sum
_{n\geq0}\en^{{-(p(\theta+1)/2-3)}} \co^{p/2} \sum
_{\tau_{i-1}^n< T} \bigl(\en^{2\theta}+\bigl(\Delta\tau
^n_i\bigr)^\theta \bigr)^{p/2} \bigl(
\Delta\tau^n_i\bigr)^{p/2-1}\Delta
\tau^n_i
\\
&\leq&\sum_{n\geq0}\en^{{-(p(\theta+1)/2-3)}}
\co^{p/2}2^{p/2-1} T \bigl(\en ^{p\theta}
D_n^{(p/2-1)}+D_n^{((\theta+1)p/2-1)} \bigr)
\\
&\leq&\co^{p/2}2^{p/2-1} T \Bigl( \Bigl(\sup_{n\geq0}
\en \Bigr)^{p\theta
/2}{\cal D}^{(p/2-1)}+{\cal D}^{((\theta+1)p/2-1)} \Bigr)<+
\infty\qquad\mbox{a.s.}
\end{eqnarray*}
Second for the case \AS--\AN, setting $D_n^{(q)}:= \sum_{\tau
_{i-1}^n< T}(\Delta\tau^n
_{i})^q$ for $q\geq0$, we have ${\cal D}^{(q)}:=\sum_{n\geq0}{\en
^{-2(q-1)+2\rhos}} D_n^{(q)} <+\infty$ a.s., thanks to
Proposition \ref{prop1}{(ii)}. Then we easily deduce (for any $p>0$)
\begin{eqnarray*}
\sum_{n\geq0}V_T^n
&\leq&\co^{p/2}2^{(p/2-1)_+} \sum_{n\geq0}{
\en^{-2(p(\theta
+1)/2-1)+2\rhos} }\\
&&\hspace*{80pt}{}\times \sum_{\tau_{i-1}^n< T} \bigl(
\en^{p\theta}\bigl(\Delta\tau^n_i
\bigr)^{p/2}+\bigl(\Delta \tau^n_i
\bigr)^{(\theta+1)p/2} \bigr)
\\
&=& \co^{p/2}2^{(p/2-1)_+} \bigl({\cal D}^{(p/2)}+{\cal
D}^{((\theta
+1)p/2)} \bigr)<+\infty\qquad\mbox{a.s.}
\end{eqnarray*}

\subsection{Proof of Proposition \texorpdfstring{\protect\ref{prop4}}{2.4}}
\label{appendix:lem6}
It is standard to check that $\cT^n$ is a sequence of increasing
stopping times; we skip details.
Let us justify that the size of $\cT^n$ is \emphh{a.s.} finite, for
any $n\geq
0$. For a given $n\geq0$, define the event $\cN^n:=\{{N^n_T}=+\infty
\}$.
For $\omega\in\cN^n$,
the infinite sequence $(\tau^n_i{(\omega)})_{i\geq0}$ converges, because
increasing and bounded by $T$. Thus, on $\cN^n\cap E_S$ with $E_S=\{
(S_t)_{0\leq t \leq T}$ continuous and $\sup_{0\leq t <T}\lmax
(H_t)<+\infty\}$, we have
\begin{eqnarray*}
0<\en&=& (S_{\tau^n_{i}} - S_{\tau^n_{i-1}})^* H_{\tau^n_{i-1}}
(S_{\tau^n_{i}} - S_{\tau^n_{i-1}}) \\
&\leq&\sup_{0\leq t <T}
\lmax(H_t) |S_{\tau^n_{i}} - S_{\tau^n
_{i-1}}|^2
\mathop{\rightarrow}_{i\rightarrow+\infty} 0,
\end{eqnarray*}
which is impossible. Thus, $\cN^n\subset E^c_S$ and $\P(\cN^n)=0$ since
$S$ is \emphh{a.s.} continuous and $\sup_{0\leq t <T}\lmax(H_t)$ is
\emphh{a.s.} finite.

Besides, we have $C_H:=\inf_{0\leq t < T}\lmin(H_t)>0$ \emphh{a.s.},
and we
immediately get
\[
\en^{-2}\supi\sup_{t\in(\tau^n_{i-1},\tau^n_i] }|\Delta S_t|^2
\leq C_H^{-1}\en^{-2} \supi\sup
_{t\in(\tau^n_{i-1},\tau^n_i] } \bigl(\Delta S_{t}^* H_{\tau
^n_{i-1}} \Delta
S_{t} \bigr) \leq C_H^{-1},
\]
which validates the assumption \AS.

Then, writing ${N^n_T}=1+\sum_{1\leq i \leq{N^n_T}-1} 1$, we point
out {(for
$n$ large enough so that $\en\leq1$)}
\begin{eqnarray*}
\en^{2\rhos} {N^n_T}&\leq&\en^{2}{N^n_T}\\
&\leq&\en^{2} + \sum_{1\leq i \leq{N^n_T}-1} \Delta
S_{\tau^n
_{i}}^* H_{\tau^n_{i-1}} \Delta S_{\tau^n_{i}}\leq
\en^2+\sum_{\tau_{i-1}^n< T}\Delta
S_{\tau^n
_{i}}^* H_{\tau^n_{i-1}} \Delta S_{\tau^n_{i}};
\end{eqnarray*}
using
moreover from Proposition \ref{prop3}, we know that under the
assumption \AS\ only,
\[
\sum_{\tau_{i-1}^n< T}\Delta S_{\tau^n_{i}}^*
H_{\tau^n_{i-1}} \Delta S_{\tau^n_{i}} \cvas \int_0^T
\tr \bigl(H_t \,\mathrm{ d} \langle S\rangle_t \bigr)<+
\infty.
\]
This validates the assumption \AN. 
%
%
\begin{remark}
The structure of hitting times of ellipsoids with size $\en$ has a
specific feature compared to general admissible strategies: the
assumption \AS\ entails the assumption \AN.
\end{remark}

\subsection{Proof of Lemma \texorpdfstring{\protect\ref{lem5}}{3.1}}
\label{subsection:proof:lem5}
We split the proof into several steps.
Let
\[
h\dvtx\cases{
\R^d\times\R_+  \rightarrow\R, \vspace*{2pt}\cr
\displaystyle(\lambda,y) \mapsto(4+d)y- \sum_{i=1}^d \sqrt{y^2+4\lambda_i^2}.}
\] Assume for a while that:
\begin{enumerate}[$(\star)$]
\item[$(\star)$]
\begin{enumerate}[(a)]
\item[(a)] for any $\lambda\in\R^d$, there exists a unique nonnegative
root $y_\lambda$ satisfying $h(\lambda,y_{\lambda})=0$;
\item[(b)]$y_0=0$; $\lambda\neq0 \Rightarrow y_\lambda>0$;
\item[(c)] the mapping $\lambda\mapsto y_\lambda$ is continuous.
\end{enumerate}
\end{enumerate}
\textit{Necessary conditions on the spectrum of $x(c)$.}
Let $\Diag$ denote the set of $d\times d$ diagonal matrices.
Take $c\in\S^d(\R)$ and let $x(c)\in\mathcal{S}^d_{+}(\R)$ be a
solution (whenever it exists) to \eqref{eq5}. Then by the spectral
theorem, $x(c)$ is diagonalizable: there exists an orthogonal matrix
$p_{x(c)}$ such that $p_{x(c)}^*x(c)p_{x(c)}\in\Diag$. Equation
\eqref
{eq5} is stable by unitary transformation
%
\begin{eqnarray}
\label{eq:matrice} && 2 \tr \bigl(p_{x(c)}^*x(c)p_{x(c)}
\bigr)p_{x(c)}^*x(c)p_{x(c)}+4 \bigl(p_{x(c)}^*x(c)p_{x(c)}
\bigr)^2
\nonumber
\\[-8pt]
\\[-8pt]
\nonumber
&&\qquad=p_{x(c)}^*c^2p_{x(c)} \in
\Diag.
\end{eqnarray}
The diagonal elements of $p_{x(c)}^*c^2p_{x(c)}$ must be the
eigenvalues of $c^2$, that is the square of the eigenvalues of $c$
[which is in $\mathcal{S}^d(\R)$]. Identifying the diagonal elements
from \eqref{eq:matrice}
gives a relation between the spectra of $c$ and $x(c)$,
\[
2\tr \bigl(x(c) \bigr) \lambda_i \bigl(x(c) \bigr) + 4
\lambda_i \bigl(x(c) \bigr)^2 = \lambda_i(c)^2,\qquad
1\leq i\leq d.
\]
Thus, the nonnegative eigenvalues of $x(c)$ must satisfy
$\lambda_i(x(c))=\break (-\tr(x(c))+\sqrt{\tr(x(c))^2+4\lambda
_i(c)^2})/4$. By
summing over $i=1,\ldots,d$, we obtain an implicit equation for $\tr
(x(c))$, which is
$h(\lambda(c),\tr(x(c)))=0$. By $(\star)$, there is a unique
solution and
%
\begin{equation}
\label{eq:tr:xc} \tr \bigl(x(c) \bigr)=y_{\lambda(c)}.
\end{equation}
Thus, we have proved that the eigenvalues of $x(c)$ must be
%
\begin{equation}
\label{eq:ev:xc} \lambda_i \bigl(x(c) \bigr)=\frac{-y_{\lambda
(c)}+\sqrt{y_{\lambda
(c)}^2+4\lambda
_i(c)^2}}{4}.
\end{equation}
\textit{Existence/uniqueness of solution to \protect\eqref{eq5}.}
Take $c\in\S^d(\R)$. Starting from \eqref{eq5}, owing to \eqref
{eq:tr:xc} $x(c)$ must solve
\[
\bigl(2 x(c)+\tfrac{ 1 }{2 }y_{\lambda(c)}I_d
\bigr)^2=\tfrac{ 1 }{4
}y^2_{\lambda(c)}I_d+c^2.
\]
The matrix $c^2+\frac{ 1 }{4 }y^2_{\lambda(c)}I_d$ is symmetric
nonnegative-definite, and thus it has a unique square-root (symmetric
nonnegative-definite matrix) \cite{HJ90}, Theorem 7.2.6, page 405, and
we obtain
%
\begin{equation}
\label{eq:sol:xc} x(c):= -\frac{y_{\lambda(c)}}4I_d+\frac{1}2
\biggl(\frac{y_{\lambda
(c)}^2}4I_d+c^2 \biggr)^{1/2}.
\end{equation}
The uniqueness is proved. It is now easy to check that $x(c)$ given in
\eqref{eq:sol:xc} solves~\eqref{eq5}, using the implicit equation
satisfied by $\tr(x(c))$. Last, $\lmin(c^2)>0$ if and only if $\lmin
(x(c))>0$ [owing to \eqref{eq:ev:xc}].

\textit{Continuity}. From Hoffman and Wielandt's theorem \cite
{HJ90}, page
368, the function $c\mapsto\lambda(c)$ is continuous on $\S
^d(\R
)$ into $\R^d$. Hence, combined with $(\star)(c)$, we obtain the
continuity of $c\mapsto y_{\lambda(c)}$ on $\S^d(\R)$ into $\R$.

Then, the continuity of $x(\cdot)$ at $c_0=0$ easily follows since as
$c\rightarrow0$, $y_{\lambda(c)}\rightarrow y_0=0$ and $\lambda
(x(c))\rightarrow0$ [using \eqref{eq:ev:xc}]: thus $x(c)\rightarrow0=x_0$.
For $c_0\neq0$, we invoke the property that $c\mapsto c^{1/2}$ is
locally lipschitz (and even analytic) on $\S^d_{++}(\R)$ into $\S
^d_{++}(\R)$ (\cite{SV97}, Lemma 5.2.1 page 131): we use this with
$\frac
{y_{\lambda(c)}^2}4I_d+c^2 \in\S^d_{++}(\R)$ for $c$ close enough to
$c_0$ (using $y_{\lambda(c)}>0$ for $c\neq0$). In view of \eqref
{eq:sol:xc}, the continuity of $x(\cdot)$ at $c_0\neq0$ follows.

\textit{Proof of $(\star)$.} $h$ is continuous on $\R^d\times
[0,\infty[$ into $\R$. Moreover:
\begin{itemize}
\item$h(\lambda,0)=-2\sum_{i=1}^d |\lambda_i|\leq0$ and $\lim_{y\rightarrow+\infty} h(\lambda,y)=+\infty$;
\item$h$ is continuously differentiable on $\R^d\times\,]0,\infty[$;
\item$D_y h(\lambda, y)=4+d - \sum_{1\leq j\leq d}\frac{y}{\sqrt {y^2+4\lambda
_i^2}}\geq4$, implying that $y\mapsto{h}(\lambda,y)$ is (strictly)
increasing.
\end{itemize}
Then, there is a unique $y_{\lambda}\in\R_+$ such that $h(\lambda,y_{\lambda})=0$. We point out at first glance,
$\lambda\neq0 \Leftrightarrow y_{\lambda}>0.$
The continuity of $y_{\cdot}$ is proved on $\R_*^d$ on the one hand, and at
$0$ on the other hand.
\begin{itemize}
\item On $\R_*^d\times\,]0,+\infty[ \dvtx D_y h(\lambda,y)$ exists
and is
nonzero: then by the implicit function theorem, there exists an open
set $U\subset\R^d_*$ containing $\lambda$ and an open set $V\subset
]0,+\infty[$ containing $y_{\lambda}$ such that $y$ is continuously
differentiable from $U$ to $V$. This proves the continuously
differentiability of $y_{\cdot}$ in $\R^d_*$.
\item At $\lambda=0 \dvtx h((|\lambda|)_{1\leq i\leq d},y)\leq
h(\lambda,y)$ and $ y \geq\frac{d|\lambda|}{\sqrt{4+2d}}\Leftrightarrow
h((|\lambda|)_{1\leq i\leq d},y)\geq0$. It implies $0\leq y_{\lambda
}\leq\frac{d|\lambda|}{\sqrt{4+2d}}$ and $\lim_{|\lambda
|\rightarrow0}
y_{\lambda}=0$.
\end{itemize}
That completes the continuity of $\lambda\mapsto y_{\lambda}$ on $\R^d$
and by the previous discussion, the proof of the lemma.\qed

\subsection{Proof of Lemma \texorpdfstring{\protect\ref{lem2}}{3.2}}
\label{proof:lem2}
We have $\langle R^n \rangle_T=\int_0^T |\sigma_t^*(D_x
u_t- D_x
u_{\varphi(t)}- D^2_{xx} u_{\varphi(t)} \Delta S_t) |^2\dt
$: to prove the result, we aim at performing a Taylor expansion using
\Au,
that is, derivatives of $u$ are \emphh{a.s.} finite in a small tube around
$(t,S_t,Y_t)_{0\leq t \leq T}$. Because of this local assumption, a
careful treatment is required, which we now detail. In view of \Au,
there exists $\Omega_\D$ such that $\P(\Omega_\D)=1$ and for every
$\omega\in\Omega_\D$ there is $\delta(\omega)>0$ such that
\[
|{\cal A} u|_{\delta}(\omega):=\sup_{0\leq t< T}\sup
_{|x-S_t(\omega
)|\leq\delta(\omega), |y-Y_t(\omega)|\leq\delta(\omega)} \bigl|\A u(t,x,y) \bigr|<+\infty
\]
for any $\A\in\D:= \{D^2_{x_j x_k}, D^3_{x_j x_k x_l}, D^2_{t
x_j}, D^2_{x_j y_m} \dvtx1\leq j,k,l \leq d,1\leq m\leq d' \}$.

Since $\supi\Delta\tau^n_i\cvas0$ and $(S_t,Y_t)_{0\leq t \leq T}$ are
\emphh{a.s.} continuous on the compact interval $[0,T]$, there exists
$\Omega_\C
$ with $\P(\Omega_\C)=1$ such that for every $\omega\in\Omega_\C$,
there is $p(\omega)\in\mathbb{N}$ such that $\forall n\geq p(\omega)$,
\[
\Bigl( \sup_{0\leq s,t\leq T, |t-s|\leq\supi\Delta\tau^n
_i}|S_t-S_s|\lor
|Y_t-Y_s| \Bigr) (\omega)\leq\delta(\omega).
\]
Hence for $\omega\in\Omega_\D\cap\Omega_\C$, let $n\geq p(\omega
)$, $i\in\{1,\ldots,{N^n_T}\}$ and $t\in[\tau^n_{i-1},\tau^n_i]$,
and write
\begin{eqnarray*}
&&D_x u(t,S_t,Y_t) - D_x u
\bigl( \tau^n_{i-1},S_{\tau^n_{i-1}},Y_{\tau^n_{i-1}}
\bigr) - D^2_{xx} u\bigl(\tau^n_{i-1},S_{\tau^n_{i-1}},Y_{\tau^n_{i-1}}
\bigr)\Delta S_t
\\
&&\qquad= \bigl[D_x u(t,S_t,Y_t)-
D_x u\bigl(\tau^n_{i-1},S_t,Y_t
\bigr) \bigr]\\
&&\qquad\quad{} + \bigl[D_x u\bigl(\tau^n_{i-1},S_t,Y_t
\bigr) - D_x u\bigl(\tau^n_{i-1},S_t,Y_{\tau^n_{i-1}}
\bigr) \bigr]
\\
&&\qquad\quad{} + \bigl[D_x u\bigl(\tau^n_{i-1},S_t,Y_{\tau^n_{i-1}}
\bigr) - D_x u\bigl(\tau^n_{i-1},S_{\tau^n_{i-1}},Y_{\tau^n_{i-1}}
\bigr)\\
&&\hspace*{83pt}\qquad\quad{}-D^2_{xx} u\bigl(\tau^n_{i-1},S_{\tau^n
_{i-1}},Y_{\tau^n_{i-1}}
\bigr)\Delta S_t \bigr].
\end{eqnarray*}
Now apply Taylor's theorem to the terms above, by observing that the
involved derivatives of $u$ are locally bounded by the (\emphh{a.s.} finite)
random variable $C_u:=\max_{\A\in\D}|{\cal A} u|_{\delta}$,
\begin{eqnarray*}
&& \bigl|D_x u(t,S_t,Y_t)- D_x u
\bigl( \tau^n_{i-1},S_{\tau^n_{i-1}},Y_{\tau^n
_{i-1}}
\bigr) - D^2_{xx} u\bigl(\tau^n_{i-1},S_{\tau^n_{i-1}},Y_{\tau^n_{i-1}}
\bigr)\Delta S_t \bigr|
\\
&&\qquad\leq\sqrt d C_u \biggl(\bigl(t-\tau^n_i
\bigr) + \sqrt{ d '} |Y_{t}-Y_{\tau^n
_{i-1}}| +
\frac{ d }{2 }|\Delta S_t|^2 \biggr).
\end{eqnarray*}
Plugging this estimate into $\langle R^n \rangle_T$
and using that $Y$ is nondecreasing, we derive that \emphh{a.s.} for
$n$ large enough,
\begin{eqnarray*}
&&{\en^{2-4\rhos}} \bigl\langle R^n \bigr\rangle_T
\\
&&\qquad\leq3 d C_u^2 \supt|\sigma_t|^2
{\en^{2-4\rhos}}\sum_{\tau_{i-1}^n< T} \biggl(\bigl(\Delta
\tau^n_i\bigr)^3 + d'|\Delta
Y_{\tau^n_i}|^2\Delta\tau^n_i
 \\
 &&\hspace*{141pt}\qquad\quad{}+ \frac{d^2}4\Delta\tau^n_i\supit|\Delta
S_t|^4 \biggr).
\end{eqnarray*}
To prove the \emphh{a.s.} convergence of the upper bound to 0, we separately
analyze each of the three contributions:
\begin{itemize}
\item
{$\en^{2-4\rhos}\sum_{\tau_{i-1}^n< T}(\Delta\tau^n_i)^3\leq\en
^{2-4\rhos} {N^n_T}\supi
(\Delta\tau^n_i)^3\leq\co\en^{4-3\rhos}\cvas0$ by Corollary \ref
{cor2}{(ii)} with $\rho=\frac{4}3 -\rhos>0$; see \AN.}
\item Combining \AY\ and Corollary \ref{cor2}{(ii)} with {$\rho
=\frac{\rho_Y}2-2(\rhos-1)>0$}, we easily obtain
\begin{eqnarray*}
{\en^{2-4\rhos}}\sum_{\tau_{i-1}^n< T}|\Delta
Y_{\tau
^n_i}|^2\Delta \tau^n_i&\leq&\sum
_{j=1}^{d'} \bigl(Y^j_{T}-Y^j_{0}
\bigr) { \en^{2-4\rhos}}\supi\bigl|\Delta Y^j_{\tau^n
_i}\bigr|\supi
\Delta \tau^n_i
\\
&\leq&\sqrt{d'} |Y_{T}-Y_{0}|
C_0 {\en^{2-4\rhos}}\en^{\rho_Y} \en^{2-\rho}\\
&\leq&
\co{\en^{\rho_Y/2-2(\rhos-1)}}\cvas0.
\end{eqnarray*}
\item Using \AS, ${\en^{2-4\rhos}}\sum_{\tau_{i-1}^n< T}\Delta
\tau^n_i\supit|\Delta
S_t|^4\leq C_0 {\en^{6-4\rhos}} T \cvas0$
{since $\rhos<\frac{3} 2$.}
\end{itemize}
All these convergences lead to the results. 

\subsection{Assumption \Au}
\label{subsection:Au}
We show that assumption \Au\ is satisfied in most usual situations,
even if the payoff $g$ is not smooth. Actually, we have not been able
to exhibit an example of $g$ for which \Au\ does not hold. The
following discussion should convince the reader that finding a
counter-example is far from being straightforward, but we conjecture
that it is possible.

\textit{Vanilla option in Black--Scholes model.} For \emph{pedagogic
reasons}, we start with the one-dimensional log-normal model $\dS
_t=\sigma S_t \dB_t$ ($\sigma>0$). Consider first the Call option with
strike $K>0$: for $t<T$ we have $D_x u(t,x)=\cN(\frac{ \log(x/K)
}{\sigma\sqrt{T-t}}+\frac{ 1 }{ 2}\sigma\sqrt{T-t} )\in[0,1]$
where $\cN(\cdot)$ is the c.d.f. of the standard Gaussian law. The second
derivative writes
\[
D^2_{xx} u(t,x)=\frac{1 }{\sigma x\sqrt{ 2\pi(T-t)} }\exp \biggl(-
\frac{ 1
}{2 } \biggl[\frac{ \log(x/K) }{\sigma\sqrt{T-t}}+\frac{ 1 }{
2}\sigma\sqrt{T-t}
\biggr]^2 \biggr);
\]
thus bounding the exponential term by $1$, we have for any given
$t_0<T$ $\lim_{\delta\rightarrow0}\sup_{0\leq t \leq
t_0}\sup_{|x-S_t|\leq\delta} |D^2_{xx} u(t,x)|\leq\frac{1 }{\sigma
\inf_{0\leq t \leq T}S_t \sqrt{ 2\pi(T-t_0) }}<+\infty$.\break  It shows
that an
\emphh{a.s.} finite bound on the second derivative is available
provided that
the time to maturity does not vanish. For the third derivative, this is
similar: indeed using $\sup_{y\in\R}e^{y^2/4}|\partial_y
(e^{-y^2/2})|=\sup_{y\in\R} |y| e^{-y^2/4}=\sqrt{2}e^{-1/2}\leq1$,
we deduce
\[
\bigl|D^3_{xxx} u(t,x)\bigr|\leq\frac{ 1+ \sigma\sqrt{ T}}{x^2\sqrt{2\pi}
\sigma^2 (T-t) } \exp \biggl(-
\frac{ 1 }{4 } \biggl[\frac{ \log(x/K) }{\sigma\sqrt{T-t}}+\frac
{ 1 }{ 2}\sigma\sqrt{T-t}
\biggr]^2 \biggr),
\]
and as before $\lim_{\delta\rightarrow0}\sup_{0\leq t \leq t_0}\sup_{|x-S_t|\leq\delta} |D^3_{xxx} u(t,x)|
<+\infty$ for any given \mbox{$t_0<T$}.

The next step consists in deriving \emphh{a.s.} upper bounds on
derivatives for
arbitrary small time to maturity. We take advantage of the property $\P
(S_T\neq K)=1$, which implies (by \emphh{a.s.} continuity of $S$) that
for $\P
$-a.e. $\omega$ there exists $t_0(\omega)\in[0,T[$ such that $\inf_{t_0(\omega)\leq t \leq T}|S_t(\omega)-K|\geq|S_T(\omega
)-K|/2:=2\delta_0(\omega)>0$.
Then, for $t\in[t_0(\omega),T]$ and $\delta\leq\delta_0\land
[2^{-1}\inf_{0\leq t \leq T}S_t]$, we have $\inf_{|x-S_t|\leq\delta
}|\log(x/K)|\geq\inf_{u>0\dvtx|u-1|\geq\delta_0/K}|\log
(u)|:=c(\omega)>0$ and\break  $\inf_{|x-S_t|\leq\delta}x\geq S_t/2$:
therefore using the inequality $-(\alpha+\beta)^2\leq-\frac{ \alpha
^2}{2 }+\beta^2$, we obtain, for $t\in[t_0(\omega),T[$
\[
\sup_{|x-S_t|\leq\delta} \bigl|D^2_{xx} u(t,x)\bigr| \leq
\frac{2 }{\sigma S_t\sqrt{ 2\pi(T-t)} }\exp \biggl(-\frac
{c^2(\omega
) }{4\sigma^2(T-t)} +\frac{1}8
\sigma^2 T \biggr).
\]
Observe that $c(\omega)>0$ implies that the above upper bound\break  converges
to 0 as $t\rightarrow T$: thus, we have completed the proof of\break  $\lim_{\delta\rightarrow0}\sup_{0\leq t< T}\sup_{|x-S_t|\leq\delta}
|D_{xx}^2 u(t,x) |<+\infty$ \emphh{a.s.} For the third derivative,
similarly
we obtain
for $t\in[t_0(\omega),T[$ and $\delta\leq\delta_0(\omega)\land
[2^{-1}\inf_{0\leq t \leq T}S_t(\omega)]$
\[
\sup_{|x-S_t|\leq\delta} \bigl|D^3_{xxx} u(t,x)\bigr|\leq
\frac{4 (1+\sigma
\sqrt T) }{S^2_t\sqrt{ 2\pi}\sigma^2(T-t) }\exp \biggl(-\frac
{c^2(\omega)
}{8\sigma^2(T-t)} +\frac{1}{16}
\sigma^2 T \biggr),
\]
and we conclude as for the second derivative. To derive the property
for $D^2_{tx} u$, we use the relation $D^2_{tx} u=-\frac{ 1 }{2
}\sigma
^2 x^2 D^3_{xxx}u-\sigma^2 x D^2_{xx}u$. Finally, \Au\ is proved for
the call option (and thus for the put option).

The same argumentation can be applied for the digital call option which
payoff is of the form $g(x)=\1_{x\geq K}$: indeed, the derivatives of
$u$ blow up only at the discontinuity point $K$ which has null
probability for the law of $S_T$. \Au\ holds for digital options.

\textit{Vanilla option in general local volatility model.} The
previous arguments are based on the explicit Black--Scholes formula for
call and digital call options, but we can generalize them to more
general models and payoffs and handle derivatives at any order. Denote
by $X^j=\log(S^j)$ ($1\leq j \leq d$) the log-asset price in a
diffusion model, and assume that $\dX_t=b^X(t,X_t)\dt+\sigma
^X(t,X_t)\dB
_t$ for coefficients $b^X$ and $\sigma^X$ of class $\mathcal
{C}_b^{\infty} ([0,T]\times\R^d )$ (bounded with bounded
derivatives). The price function in the log-variables is then
$v(t,x):=u(t,\exp(x^1),\ldots,\exp(x^d))=\E(g(S_T)|S^j_t=\exp
(x^j),1\leq
j \leq d):=\break \E(G(X_T)|X_t=x)$. We first consider the simple case of
${\mathcal C}^\infty$-payoff $G$ with exponentially bounded
derivatives: for any $k\geq0$, there is a constant $C^G_k\geq0$ such
that $|D^k_x G (x)|\leq C^G_k \exp(C^G_k|x|)$ for $x\in\R^d$. In this
case, a direct differentiation of $\E(G(X_T)|X_t=x)$ using the smooth
flow $x\mapsto X_T^{t,x}$ \cite{kuni:84} shows the differentiability of
$v$ w.r.t. the space variable with derivatives bounded on compact
subsets of $[0,T]\times\R^d$; in addition the time smoothness is
obtained using It\^{o}'s formula; these arguments are standard and we
skip details. \Au\ is proved for these smooth payoffs.

Now we tackle the case of discontinuous payoffs of the form $G(x)=\1
_{x\in\D}\varphi(x)$ for a closed set $\D\subset\R^d$ and a
$\mathcal
{C}^{\infty}$-function $\varphi$ with exponentially bounded
derivatives: observe that by combining the analysis for smooth payoffs
and that for discontinuous ones will allow to cover a quite large class
of $g$ satisfying \Au\ (such as call/put, digital call/put, exchange
call, digital exchange call and so on). We assume that a uniform
ellipticity assumption is satisfied:
$\inf_{0 \leq t \leq T,x\in\R^d}\inf_{|\xi|=1}\xi\cdot[\sigma
^X(\sigma
^X)^*](t,x)\xi>0$. In this setting, $v(t,x)=\int_{\R^d} \1_{z\in\D
}p(t,x,T,z)\varphi(z) \,\dz$ where $p$ is the transition density function
of $X$, which is smooth and satisfies to Aronson-type estimates (\cite{frie:64},
Theorem 8, page 263): for any $i\geq0$ and any
differentiation index $\alpha$, there exists a constant $C_{i,\alpha
}=C_{i,\alpha}(T,b^X,\sigma^X)>0$ such that
\[
\bigl|D^{i,\alpha}_{tx} p(t,x,T,z)\bigr|\leq C_{i,\alpha}(T-t)^{-(d+2i+|\alpha
|)/2}
\exp \bigl(-|x-z|^2/ \bigl[C_{i,\alpha}(T-t) \bigr] \bigr)
\]
for any $0\leq t <T$, $x
\in\R^d$, $z\in\R^d$. From the integral representation of $v$, it
readily follows that
\begin{eqnarray*}
&&\bigl|D^{i,\alpha}_{tx} v(t,x)\bigr|\\
&&\qquad\leq C_{i,\alpha}(T-t)^{-(2i+|\alpha
|)/2}
\int_{\R^d}C^\varphi_0 e^{C_0^\varphi|z|}
(T-t)^{-d/2}e^{-|x-z|^2/[C_{i,\alpha}(T-t)]}\,\dz
\\
&&\qquad\leq C_{i,\alpha}(T-t)^{-(2i+|\alpha|)/2}C^\varphi_0
e^{C_0^\varphi|x|} \int_{\R^d}e^{C_0^\varphi\sqrt
{T}|w|}e^{-|w|^2/C_{i,\alpha}}
\,\mathrm{ d}w,
\end{eqnarray*}
which proves locally uniform bounds on derivatives provided that the
time to maturity remains bounded away from 0. To handle the case
$t\rightarrow T$, we additionally assume that \emph{the boundary $\dD$
of $\D$ is Lebesgue-negligible} (thus including usual situations but
excluding Cantor like sets; see \cite{dibe:02}, page 114): thus for $\P
$-a.e. $\omega$, the distance to the boundary (a closed set) is
positive, that is, $\delta_0(\omega):=\frac{ 1 }{4 }d(X_T(\omega
),\dD
)>0$, and there exists $t_0(\omega)\in[0,T[$ such that $\inf_{t_0(\omega
)\leq t \leq T}d(X_t(\omega),\dD)\geq3\delta_0(\omega)$ [we recall
that the distance function $x\mapsto d(x,\dD)$ is Lipschitz
continuous]. Now, let $\omega$ be given as above; by the smooth version
of the Urysohn lemma~\cite{dieu:90:1}, page~90, Chapter IV, there exists a smooth
function $\xi$ (depending on $\omega$) such that $\1_{x \in\D,
\delta
_0\leq d(x,\dD)}\leq\xi(x)\leq\1_{x\in\D}$. Decompose the price
function into two parts $v=v_1+v_2$ with
\begin{eqnarray*}
v_1(t,x)&:=&\int_{\R^d} \1_{z\in\D}p(t,x,T,z)
\varphi(z) \xi(z)\,\dz,\\
  v_2(t,x)&=&\int_{\D}
p(t,x,T,z)\varphi(z) \bigl(1-\xi(z) \bigr) \,\dz.
\end{eqnarray*}
We easily handle the derivatives of $v_1$ using the first case of
smooth functions since $\1_{\D}\varphi\xi=\varphi\xi\in{\mathcal
C}^\infty$ with exponentially bounded derivatives. Regarding $v_2$,
observe that we integrate over the $z$ such that $z\in\D$ and
$d(z,\,\dD
)<\delta_0$; for such $z$, for $t\in[t_0,T[$ and $|x-X_t|\leq\delta
\leq
\delta_0$,
we have $|x-z|\geq d(X_t,\dD)-|x-X_t|-d(z,\dD)\geq\delta_0$ and thus
\begin{eqnarray*}
&&\sup_{|x-X_t|\leq\delta} \bigl|D^{i,\alpha}_{tx}v_2(t,x)\bigr|
\\
&&\qquad\leq\sup_{|x-X_t|\leq\delta} \int_\D
C^\varphi_0 e^{C_0^\varphi|z|} C_{i,\alpha}(T-t)^{-(d+2i+|\alpha
|)/2}e^{-|x-z|^2/[2C_{i,\alpha
}(T-t)]}\\
&&\hspace*{50pt}\quad\qquad{}\times e^{-\delta_0^2/[2C_{i,\alpha}(T-t)]}
\,\dz
\\
&&\qquad\leq C_{i,\alpha}(T-t)^{-(2i+|\alpha|)/2} e^{-\delta
_0^2/[2C_{i,\alpha}(T-t)]}
C^\varphi_0 e^{C_0^\varphi(|X_t|+\delta_0)}\\
&&\qquad\quad{}\times  \int_{\R
^d}e^{C_0^\varphi\sqrt{T}|w|}e^{-|w|^2/[2C_{i,\alpha
}]}
\,\mathrm{ d}w.
\end{eqnarray*}
The above upper bound converges to 0 as $t\rightarrow T$, and the proof
of \Au\ is complete.

Interestingly, we can weaken the ellipticity assumption into a
hypoellipticity assumption: indeed, our analysis essentially relies on
transition density estimates in small time and away from the diagonal.
These estimates are available in the hypoelliptic homogeneous diffusion
case (\cite{kusu:stro:85}, Corollary 3.25) and in the inhomogeneous case
\cite{catt:mesn:02}, Assumption (1.10).

\textit{Asian option in general local volatility model.} The payoff
is of the form $g(S_T,I_T)$ where $I_T=\int_0^T S_t \dt$ and $S$ is a
one-dimensional homogeneous diffusion $\dS_t=\sigma
(S_t)\dB_t$.
The analysis is reduced to the previous case of vanilla option by
considering the 2-dimensional diffusion $(S_t,I_t)_{0\leq t \leq T}$:
it is not elliptic but hypoelliptic \cite{kusu:stro:85} provided that
$\sigma$ is smooth and that $\sigma(x)>0$ for $x\in I$ where
$I\subset
\R$ is given by $\P(\forall t\in[0,T]\dvtx X_t\in I)=1$ (in usual cases,
$I=\,]0,+\infty[$). It includes the Black--Scholes model and any model
with local volatility bounded away from 0 and smooth. We skip details.

\textit{Lookback option in Black--Scholes model}. The payoff is of
the form $S_T-m\land\min_{0\leq t \leq T}S_t$ or $M\lor\max_{0\leq t
\leq T}S_t-S_T$ for lookback call or put, $(M\lor\max_{0\leq t \leq
T}S_t-K)_+$ or $(K-m\land\min_{0\leq t \leq T}S_t)_+$ for call on
maximum or on minimum, $(S_T-\lambda m\land\min_{0\leq t \leq
T}S_t)_+$ (with $\lambda>1$) or $(\lambda M\lor\max_{0\leq t \leq
T}S_t-S_T)_+$ (with $\lambda<1$) for partial lookback call or put. In
all these cases, Black--Scholes-type formulas are available in closed
forms \cite{conz:visw:91}. Then it is straightforward to check that
\Au\ is satisfied, and this is essentially based on the property that
under the assumption of nonzero volatility, the joint law $(S_T,\max_{0\leq t \leq T}S_t, \min_{0\leq t \leq T}S_t)$ has a density (derived
from \cite{RY99}, Exercise 3.15), implying that the events on which the
derivatives may blow up (such as $\{
S_T=\min_{0\leq t \leq T}S_t\}, \ldots$) have zero probability.
\end{appendix}

\section*{Acknowledgments}
The author is grateful to the Association Nationale
de la Recherche Technique and GDF SUEZ for their financial support.

%

%


\printaddresses

\end{document}